\documentclass{amsart}
\usepackage[all]{xy}
\xyoption{arc}

\input epsf.tex
\newcommand{\evects}{$$
\xymatrix{
 & z^kv_1 \ar@{-}[d]|{\textstyle d} & z^kv_3 \ar@{-}[d]|{ \textstyle c} & \\
z^{-j}v_2 \ar@{-}[r]|{ \textstyle a} & v_4 \ar@{-}[r]|{ \textstyle b} \ar@{-}[d]|{ \textstyle c} & v_2 \ar@{-}[d]|{ \textstyle d} &
\ar@{-}[l]|{ \textstyle a} z^jv_4 \\ 
z^{-j}v_3 \ar@{-}[r]|{ \textstyle b} & v_1 \ar@{-}[r]|{ \textstyle a} & v_3 & \ar@{-}[l]|{ \textstyle b} z^jv_1 \\
 & z^{-k}v_4 \ar@{-}[u]|{ \textstyle d} & z^{-k}v_2 \ar@{-}[u]|{ \textstyle c} & \\
}
$$
}

\newcommand{\ellipsi}{\begin{picture}(0,0)%
\epsfbox{ell.pstex}%
\end{picture}%
\setlength{\unitlength}{0.00083300in}%
\begingroup\makeatletter\ifx\SetFigFont\undefined
\def\x#1#2#3#4#5#6#7\relax{\def\x{#1#2#3#4#5#6}}%
\expandafter\x\fmtname xxxxxx\relax \def\y{splain}%
\ifx\x\y   
\gdef\SetFigFont#1#2#3{%
  \ifnum #1<17\tiny\else \ifnum #1<20\small\else
  \ifnum #1<24\normalsize\else \ifnum #1<29\large\else
  \ifnum #1<34\Large\else \ifnum #1<41\LARGE\else
     \huge\fi\fi\fi\fi\fi\fi
  \csname #3\endcsname}%
\else
\gdef\SetFigFont#1#2#3{\begingroup
  \count@#1\relax \ifnum 25<\count@\count@25\fi
  \def\x{\endgroup\@setsize\SetFigFont{#2pt}}%
  \expandafter\x
    \csname \romannumeral\the\count@ pt\expandafter\endcsname
    \csname @\romannumeral\the\count@ pt\endcsname
  \csname #3\endcsname}%
\fi
\fi\endgroup
\begin{picture}(6058,2355)(1172,-2422)
\put(2476,-1261){\makebox(0,0)[lb]{\smash{\SetFigFont{10}{10}{rm}$-(a-b)^2$}}}
\put(1276,-1261){\makebox(0,0)[lb]{\smash{\SetFigFont{10}{10}{rm}$-(c+d)^2$}}}
\put(3451,-1261){\makebox(0,0)[lb]{\smash{\SetFigFont{10}{10}{rm}$0$}}}
\put(3076,-2386){\makebox(0,0)[lb]{\smash{\SetFigFont{10}{10}{rm}$(\theta_0,\phi_0)$}}}
\put(3001,-211){\makebox(0,0)[lb]{\smash{\SetFigFont{10}{10}{rm}$(-\theta_0,-\phi_0)$}}}
\put(4201,-1261){\makebox(0,0)[lb]{\smash{\SetFigFont{10}{10}{rm}$(c-d)^2$}}}
\put(6526,-1261){\makebox(0,0)[lb]{\smash{\SetFigFont{10}{10}{rm}$(a+b)^2$}}}
\put(6601,-2386){\makebox(0,0)[lb]{\smash{\SetFigFont{10}{10}{rm}$\theta$}}}
\put(1426,-406){\makebox(0,0)[lb]{\smash{\SetFigFont{10}{10}{rm}$\phi$}}}
\put(6251,-211){\makebox(0,0)[lb]{\smash{\SetFigFont{10}{10}{rm}$E_1$}}}
\put(1576,-2161){\makebox(0,0)[lb]{\smash{\SetFigFont{10}{10}{rm}$E_2$}}}
\end{picture}
}
\newcommand{\missellips}{\begin{picture}(0,0)%
\epsfbox{ell2.pstex}%
\end{picture}%
\setlength{\unitlength}{0.00083300in}%
\begingroup\makeatletter\ifx\SetFigFont\undefined
\def\x#1#2#3#4#5#6#7\relax{\def\x{#1#2#3#4#5#6}}%
\expandafter\x\fmtname xxxxxx\relax \def\y{splain}%
\ifx\x\y   
\gdef\SetFigFont#1#2#3{%
  \ifnum #1<17\tiny\else \ifnum #1<20\small\else
  \ifnum #1<24\normalsize\else \ifnum #1<29\large\else
  \ifnum #1<34\Large\else \ifnum #1<41\LARGE\else
     \huge\fi\fi\fi\fi\fi\fi
  \csname #3\endcsname}%
\else
\gdef\SetFigFont#1#2#3{\begingroup
  \count@#1\relax \ifnum 25<\count@\count@25\fi
  \def\x{\endgroup\@setsize\SetFigFont{#2pt}}%
  \expandafter\x
    \csname \romannumeral\the\count@ pt\expandafter\endcsname
    \csname @\romannumeral\the\count@ pt\endcsname
  \csname #3\endcsname}%
\fi
\fi\endgroup
\begin{picture}(4612,3314)(1922,-2919)
\put(3826,-886){\makebox(0,0)[lb]{\smash{\SetFigFont{10}{10}{rm}$E_2$}}}
\put(5551,-436){\makebox(0,0)[lb]{\smash{\SetFigFont{10}{10}{rm}$E_1$}}}
\put(3451,-1186){\makebox(0,0)[lb]{\smash{\SetFigFont{10}{10}{rm}$0$}}}
\put(2701,-361){\makebox(0,0)[lb]{\smash{\SetFigFont{10}{10}{rm}$\phi$}}}
\put(6226,-2761){\makebox(0,0)[lb]{\smash{\SetFigFont{10}{10}{rm}$\theta$}}}
\end{picture}
}

\theoremstyle{plain}
\newtheorem{thm}{Theorem}
\newtheorem{cor}[thm]{Corollary}
\newtheorem{lem}[thm]{Lemma}
\newtheorem{prop}[thm]{Proposition}
\newtheorem{conj}{Conjecture}
\newtheorem{prob}{Open Problem}
\numberwithin{thm}{section}
\numberwithin{conj}{section}
\numberwithin{prob}{section}

\numberwithin{equation}{section}

\newcommand{\R}{{\mathbb R}}
\newcommand{\Z}{{\mathbb Z}}
\newcommand{\C}{{\mathbb C}}
\newcommand{\N}{{\mathbb N}}
\newcommand{\E}{{\bf E}}

\newcommand{\p}{\partial}

\newcommand{\length}{\textup{length}}
\newcommand{\area}{\textup{area}}
\newcommand{\Asymp}{\textup{AH}}
\newcommand{\Ent}{\textup{Ent}}
\newcommand{\ent}{\textup{ent}}
\newcommand{\av}{\textup{av}}
\newcommand{\Prob}{\textup{Prob}}

\newcommand{\LocEnt}[1]{\ent
\left(\frac{\partial #1}{\partial x},
\frac{\partial #1}{\partial y}\right)}

\newcommand{\hfsup}{\vee}
\newcommand{\hfinf}{\wedge}

\title{A Variational Principle for Domino Tilings}

\author{Henry Cohn}
\address{Department of Mathematics\\
Harvard University\\
Cambridge, MA 02138\\
USA}
\curraddr{Microsoft Research\\
One Microsoft Way\\
Redmond, WA 98052-6399}
\email{cohn@math.harvard.edu}
\author{Richard Kenyon}
\address{CNRS UMR 8628\\
Laboratoire de Topologie\\
B\^atiment 425\\
Universit\'e  Paris-11\\
91405 ORSAY\\
France}
\email{kenyon@topo.math.u-psud.fr}
\author{James Propp}
\address{Department of Mathematics\\
University of Wisconsin\\
Madison, WI 53706\\
USA}
\email{propp@math.wisc.edu}

\thanks{Cohn was supported by an NSF Graduate Research
Fellowship.  Propp was supported  by NSA grant MDA904-92-H-3060 and
NSF grant DMS92-06374, and by a grant from the MIT Class of 1922.}


\dedicatory{Dedicated to Pieter Willem Kasteleyn (1924--1996)}

\date{August 21, 2001;  this version of the paper is updated from
the published version (Journal of the AMS {\bf 14} (2001), 297--346)
by correcting typos in \eqref{fx}, \eqref{fy}, \eqref{Peqns}, \eqref{Epq},
and Subsection~\ref{entasfnofslope}, as well as a few other places that
restate these results.  The basic problem was that in \eqref{fx} and
\eqref{fy}, $s$ and $t$ were switched and there was a sign error.
This error propagated through the paper, but is easily corrected.  None
of the underlying logic needs to be changed.}

\begin{document}

\begin{abstract}
We formulate and prove a variational principle (in the sense of
thermodynamics) for random domino tilings, or equivalently for the
dimer model on a square grid.  This principle states that a typical
tiling of an arbitrary finite region can be described by a function
that maximizes an entropy integral.  We associate an entropy to every
sort of local behavior domino tilings can exhibit, and prove that
almost all tilings lie within $\varepsilon$ (for an appropriate
metric) of the unique entropy-maximizing solution.  This gives a
solution to the dimer problem with fully general boundary conditions,
thereby resolving an issue first raised by Kasteleyn.  Our methods
also apply to dimer models on other grids and their associated tiling
models, such as tilings of the plane by three orientations of unit
lozenges.
\end{abstract}

\maketitle

\begin{quote}
\textit{The effect of boundary conditions is, however, not entirely
trivial and will be discussed in more detail in a subsequent paper.}
\begin{flushright}
P.~W.~Kasteleyn, 1961
\end{flushright}
\end{quote}

\section{Introduction}
\label{introsec}

\subsection{Description of results}
\label{statementsec}

\begin{figure}
\begin{center}
\leavevmode
\epsfbox[0 0 204 204]{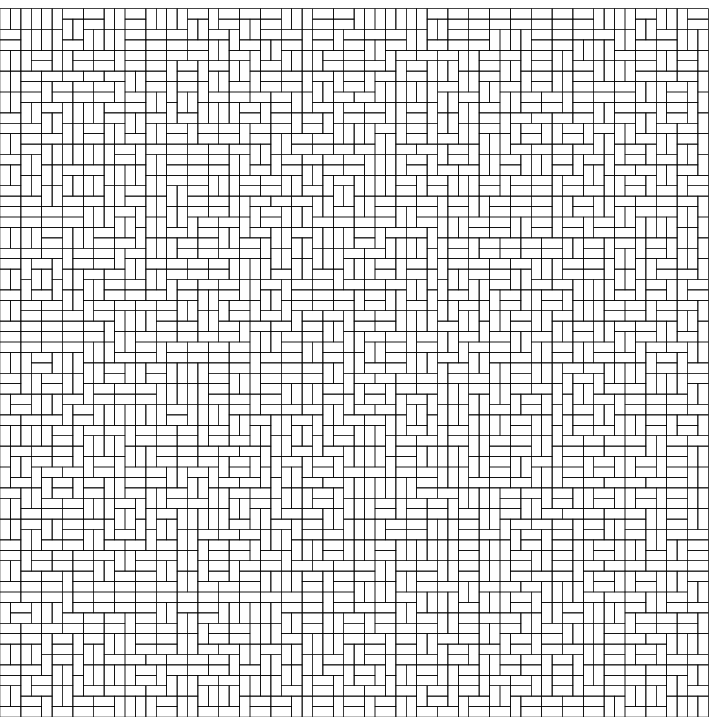}
\end{center}
\caption{A random domino tiling of a square.}
\label{fig-square}
\end{figure}

\begin{figure}
\begin{center}
\leavevmode
\epsfbox[0 0 288 288]{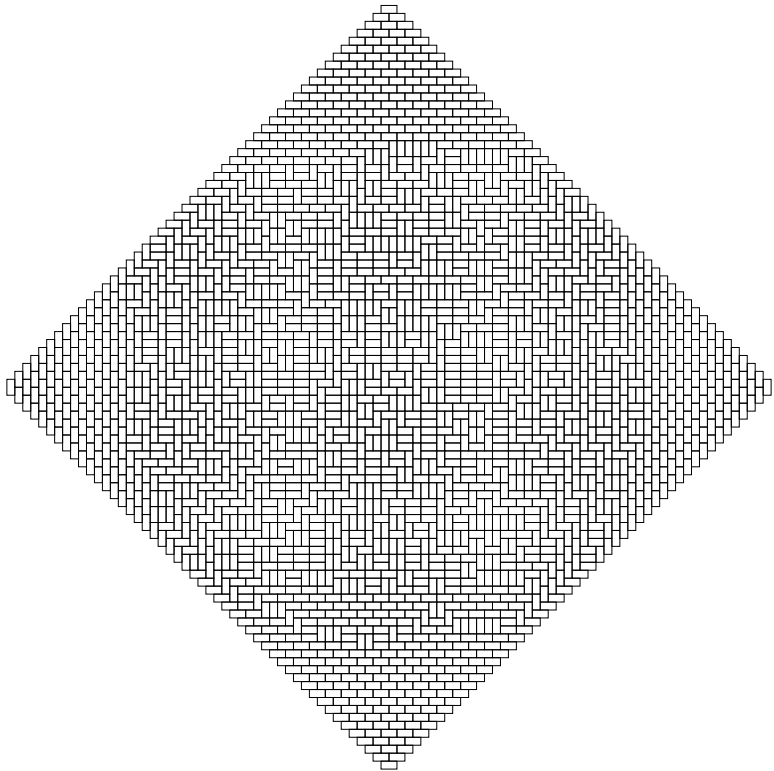}
\end{center}
\caption{A random domino tiling of an Aztec diamond.}
\label{fig-aztec}
\end{figure}

A {\it domino\/} is a $1 \times 2$ (or $2 \times 1$) rectangle,
and a {\it tiling\/} of a region by dominos
is a way of covering that region with dominos
so that there are no gaps or overlaps.
In 1961, Kasteleyn \cite{Kast1} found a formula for the number
of domino tilings of an $m \times n$ rectangle (with $mn$ even),
as shown in Figure~\ref{fig-square} for $m=n=68$.
Temperley and Fisher \cite{TF} used a different method
and arrived at the same result at almost exactly the same time.
Both lines of calculation
showed that the logarithm of number of tilings,
divided by the number of dominos in a tiling
(that is, $mn/2$), converges to
$2G/\pi \approx 0.58$ (here $G$ is Catalan's constant).
On the other hand,
in 1992 Elkies et al.\ \cite{EKLP} studied domino tilings
of regions they called Aztec diamonds
(Figure~\ref{fig-aztec} shows an Aztec diamond of order 48),
and showed that the logarithm
of the number of tilings, divided by the number of dominos, converges
to the smaller number $(\log 2)/2 \approx 0.35$.
Thus, even though the region in Figure~\ref{fig-square}
has slightly smaller area than the region in Figure~\ref{fig-aztec},
the former has far more domino tilings.
For regions with other shapes,
neither of these asymptotic formulas may apply.

In the present paper we consider
simply-connected regions of arbitrary shape.
We give an exact formula for the limiting
value of the logarithm of the number of tilings
per unit area, as a function of the shape of the boundary of the
region, as the size of the region goes to infinity.
In particular, we show that computation of this limit is
intimately linked with an understanding of long-range variations in
the local statistics of random domino tilings.
Such variations can be seen by comparing
Figures~\ref{fig-square} and~\ref{fig-aztec}.
Each of the two tilings is random in the sense that
the algorithm \cite{PW} that was used to create it generates
each of the possible tilings of the region being tiled with the same
probability.
Hence one can expect each tiling to be qualitatively typical of
the overwhelming majority of tilings of the region in question,
as is in fact the case.
Figure~\ref{fig-square} looks more or less homogeneous, but
even cursory examination of Figure~\ref{fig-aztec} shows that the
tiling manifests different gross behavior in different parts of the
region.  In particular, the tiling degenerates into a
non-random-looking brickwork pattern near the four corners of the
region, whereas closer to the middle one sees a mixture of horizontal
and vertical dominos
of the sort seen everywhere in Figure~\ref{fig-square}
(except very close to the boundary).

In earlier work \cite{CEP,JPS}
two of us, together with other researchers, analyzed
random domino tilings of Aztec diamonds in great detail, and showed
how some of the properties of Figure~\ref{fig-aztec} could be
explained
and quantified.
It was proved that
the boundary between the four brickwork
regions and the central mixed region for a randomly tiled Aztec
diamond tends in probability towards the inscribed circle
(the so-called ``arctic circle'') as the size of the diamonds
becomes large \cite{JPS}, and that even inside the
inscribed circle, the first-order local statistics
(that is, the probabilities of finding
individual dominos in particular locations)
fail to exhibit homogeneity on a macroscopic scale \cite{CEP}.

Unfortunately, the techniques of \cite{JPS} and \cite{CEP}
do not apply to general regions.  A few cases besides Aztec
diamonds have been analyzed; for example, random domino tilings
of square regions have been analyzed and
do turn out to be homogeneous \cite{BP}.
(That is, if one looks at two patches
at distance at least $d$ from the boundary of an $n \times n$ square,
with $n$ even,
the local statistics on the two patches
become more and more alike as $n$ goes to infinity,
as long as $d$ goes to infinity with $n$.)
However, before our research was undertaken no general analysis was known.

In this paper, we will demonstrate that
the behavior of random tilings of large regions
is determined by a variational (or entropy maximization) principle,
as was conjectured in Section~8 of \cite{CEP}.
We show that
{\it the logarithm of the
number of tilings, divided by the number of dominos in a tiling of $R$,
is asymptotic (when $\area(R) \rightarrow \infty$) to}
\begin{equation}
\label{supoverh}
\sup_h \iint_{\!\!R^*} \ent
\left( \frac{\partial h}{\partial x},
\frac{\partial h}{\partial y} \right) \, dx \: dy.
\end{equation}
Here the domain of integration $R^*$
is a normalized version of $R$, the function
$h$ ranges over a certain compact set of Lipschitz functions
from $R^*$ to $\R$,
and
\begin{equation}
\label{entdef}
\ent(s,t) =
\frac{1}{\pi}(L(\pi p_a) + L(\pi p_b) + L(\pi p_c) + L(\pi p_d)),
\end{equation}
where $L(\cdot)$ is the Lobachevsky function (see \cite{Mil}), defined by
\begin{equation}
\label{lobdef}
L(z) = -\int_{0}^{z} \log |2 \sin t| \: dt,
\end{equation}
and the quantities $p_a,p_b,p_c,p_d$
are determined by the equations
\begin{eqnarray}\label{s}
2(p_a-p_b)&=&t\label{fx},\\
2(p_d-p_c)&=&s\label{fy},\\
p_a+p_b+p_c+p_d&=&1\label{sump},\\
\sin(\pi p_a)\sin(\pi p_b)&=&\sin(\pi p_c)\sin(\pi p_d).\label{sinsin}
\end{eqnarray}
The quantities $p_a$, $p_b$, $p_c$, $p_d$
can be understood in terms of properties of random domino tilings
on the torus (see Subsection~\ref{interpretsec}); the quantities also
have attractively direct but still conjectural interpretations
in terms of the local statistics for random tilings of the
original planar region (see Conjecture~\ref{probconj}).

There exists a unique function $f$ that achieves the maximum in
\eqref{supoverh}.
Its partial derivatives encode
information about the local statistics
exhibited by random
tilings of the region; for example, the places $(x,y)$ where the
``tilt'' $({\partial f}/{\partial x}, {\partial f}/{\partial
y})$ is $(2,0)$, $(-2,0)$, $(0,2)$, or $(0,-2)$ correspond to the
places in the tiling where the probability of seeing brickwork
patterns goes to 1 in the limit, in a suitable sense.
The function $f$ need not be
$C^\infty$, and in fact often is not.
For example, in the case of Aztec diamonds $f$
is smooth everywhere except on the arctic circle, where it is
only $C^1$ (except at the midpoints of the sides, where it is only $C^0$).
Inside this circle, $f$ takes on a variety of values,
corresponding to the fact that different local statistics
are manifested at different locations.
In contrast, for square regions the function $f$ is constant,
corresponding to the fact that throughout the region
(except very close to the boundary),
the local statistics are constant.
(See Figures~\ref{fig-square} and~\ref{fig-aztec}.)

In general, the locations where $f$ is not smooth
correspond to the existence of a phase transition
in the two-dimensional dimer model,
closely related to a phase transition first noticed by Kasteleyn \cite{Kast2}.

The function $f$ is related to combinatorial representations
of states of the dimer model called {\it height functions\/}.
Each of the different tilings of a fixed finite simply-connected region in
the square grid has a height function which is a function from the
vertices of the region to $\Z$ satisfying certain conditions (see
Subsection~\ref{hfsec}); there is a one-to-one correspondence between
tilings and height functions (as long as we fix the height at one
point, since the actual correspondence is between tilings and height
functions modulo additive constants).
All of the height functions agree on
the boundary of the region but have different values on the interior.
In an earlier paper \cite{CEP}, it was shown that these height
functions satisfy a law of large numbers, in the sense that for each
vertex $v$ within a very large region, the height of $v$ in a random
tiling of the region is a random variable whose standard deviation is
negligible compared to the size of the region.  This suggested the
existence of a limit law, but did not indicate what the limit law was.
We show (see Theorem~\ref{maintheorem}
at the end of Subsection~\ref{sketchsec}):
{\it for every $\varepsilon>0$,
the height function of a random tiling,
when rescaled by the dimensions of the region
being tiled, is, with probability tending to 1, within $\varepsilon$ of
$f$\/}, where $f$ is the function that maximizes the double integral
in \eqref{supoverh}.
We therefore call $f$ an {\it asymptotic height function\/}.
One may qualitatively summarize our main results by
saying that the pattern governing local statistical behavior of
uniform random tilings of a region is, in the limit, the unique
pattern that maximizes the integral of the ``local entropy'' over the region.
Moreover, the value of this maximum is an asymptotic expression for
the ``global entropy'' of the ensemble of tilings of the region.

Our main theorem
(Theorem~\ref{maintheorem})
thus has two intertwined components: a law of large
numbers that describes the local statistics of random tilings of a
large region by way of an asymptotic height function $f$, and an
entropy estimate that tells us that the total number of tilings is
determined by a certain functional of that self-same function $f$.
Furthermore $f$ is precisely the unique Lipschitz function that maximizes
that functional.

We supplement our main theorem by two supporting results:
a large deviations estimate (Theorem~\ref{main})
and a PDE that the entropy-maximizing function $f$ must satisfy
(Theorem~\ref{pdethm}).
To be technically correct, we should not assert that the function $f$
satisfies the PDE everywhere, but only where the partial derivatives
are continuous and the tilt $(s,t)$ satisfies $|s|+|t| < 2$.
This proviso is necessary because singularities in
$f$ (corresponding to domain boundaries like the arctic circle) are an
essential feature of the phenomena we are studying;
in fact, we shall
find (see Section~\ref{partitionsection} and Theorem~\ref{pa})
that these singularities are related to a phase
transition manifested by the dimer model as an external field is
permitted to vary, and that domain boundaries can be viewed as a
spatial expression of that phase transition.
Our results also imply that if instead of studying the uniform
distribution on the set of {\it all\/} the tilings of a large region
$R$, one restricts one's attention to those tilings whose height
functions (suitably normalized) approximate some asymptotic height
function $h$ (which need not be the $h$ that maximizes the double
integral \eqref{supoverh}), then the entropy of this restricted
ensemble (that is, the normalized logarithm of the number of tilings)
approximates the value of the double integral.
Note that this more general
result implies that the total number of (unrestricted) tilings is at least
as large as the supremum \eqref{supoverh}.  As for integrand itself, we remark
here that $\ent(s,t)$ achieves its maximum at $(0,0)$ and that it goes to
zero as $(s,t)$ goes to the boundary curve $|s|+|t|=2$; that is, there are
many ways to tile a patch so that its average tilt is near zero, but fewer
ways to tile it as the tilt gets larger, and no ways to tile it at all
if the desired tilt $(s,t)$ fails to satisfy $|s|+|t| \leq 2$.

\subsection{Interpretation}
\label{interpretsec}

A {\it dimer configuration\/}, or {\it perfect matching\/},
of a graph is a set of edges in the graph
such that each vertex belongs to exactly one of the selected edges.
To see the equivalence between tilings by dominos and dimer
configurations on a grid (the graph-theoretic dual to the graph of
edges between lattice squares),
one need only replace each $1 \times 1$ square by a vertex and
each domino by an edge.

To explain the significance of the quantities $p_a,p_b,p_c,p_d$,
we need to digress and discuss the dimer model on a torus.
Here the graph is just like the $m \times n$ rectangular grid,
but with $m$ extra bonds that connect vertices on the left and right
and $n$ extra bonds that connect vertices on the top and bottom.
Kasteleyn
showed \cite{Kast1} that the number of dimer configurations
on this graph is governed by the same asymptotic formula as for
the dimer model on a rectangle.
In the same article
Kasteleyn considered a generalization to non-uniform distributions
on the set of tilings, obtained by assigning different ``weights''
to horizontal and vertical bonds, and letting the probability of
a particular dimer configuration be proportional to the product
of the weights of its bonds.  In the statistical mechanics
literature, such modifications of a model are sometimes conceived
of as resulting from the imposition of an external field.

Here we go one step further, and impose a field that discriminates
among four different kinds of bonds: $a$-bonds and $b$-bonds
(both horizontal) and $c$-bonds and $d$-bonds (both vertical),
staggered as in Figure~\ref{evects} from Section~\ref{detsection}
(this is different from a 4-parameter external field
that was considered by Kasteleyn \cite{Kast2}).
The field depends on four parameters $a,b,c,d$,
which are the weights associated with the respective kinds of bonds;
the probability of a dimer configuration on the torus graph
is proportional to the product of the weights of its bonds.
Taking the graph-theoretic dual, we get a non-uniform distribution
on domino tilings on the torus,
having less symmetry than the torus itself.
One way to motivate the consideration of such asymmetrical measures
is to observe that in regions like the one shown in Figure~\ref{fig-aztec},
the boundary conditions break the symmetry of the underlying square grid
in precisely this fashion,
yielding subregions in which bonds (or dominos)
of one of the four types predominate.

To avoid unnecessary complexity,
we limit ourselves to square tori ($m=n$),
but the situation for
more general
tori is much the same.

The quantities $p_a,p_b,p_c,p_d$ have a direct (and rigorously proved)
interpretation in terms of the dimer model on the $n \times n$ torus.
Suppose we are given positive real weights $a,b,c,d$
as described above.
Suppose that each of $a,b,c,d$ is less than the sum of the others.
Then there is a unique Euclidean quadrilateral of edge lengths $a,c,b,d$
(in that order) which is cyclic, that is, can be inscribed in a circle.
Define $p_a$ to be $1/(2\pi)$ times the angle
of arc of the circumscribed circle cut off by the edge $a$ of this
quadrilateral.
Similarly define $p_b,p_c$ and $p_d$.
Then we shall see (see Theorem~\ref{pa}) that $p_a$ is (in the large-$n$ limit)
the probability that a given $a$-bond
belongs to a randomly-chosen dimer configuration
(under the probability distribution determined by the weights $a,b,c,d$),
and likewise for $p_b,p_c,p_d$.
Technically, we only prove convergence for $n$ in a large subset of $\N$,
but we believe that convergence holds for all $n$.
If on the other hand $a\geq b+c+d$, then
we shall see that as $n\to\infty$, $p_a$ tends to $1$ (and $p_b,p_c,p_d$
tend to zero). A similar phenomenon occurs when $b$, $c$, or $d$ is greater
than or equal to the sum of the others.

Moreover, the quantities $s = 2(p_d-p_c)$ and $t = 2(p_a-p_b)$
have a height-function interpretation.
Since the torus can be viewed as a rolled-up plane,
every dimer configuration on the torus ``unrolls'' to give
a dimer configuration on the plane,
which (in the guise of a domino tiling of the plane)
gives rise to a height function.
The height is not in general a periodic function,
but rather increases by some amounts $H_x$ and $H_y$
as one moves $n$ vertices in the $+x$ direction
or $n$ vertices in the $+y$ direction.
Here $H_x$ and $H_y$ depend on the tiling chosen,
and thus are random variables;
the expected values of $H_x/n$ and $H_y/n$ are $s$ and $t$, respectively.

The relationship between the uniform measure on tilings
of finite simply-connected regions
and the $a,b,c,d$-weighted measure on tilings of tori
may become clearer after one has verified that
the $a,b,c,d$-weighted measure on the tilings of a
finite simply connected region, defined in the obvious fashion,
is nothing other than the uniform distribution,
as long as $ab=cd$.
To see this, one may make use of the fact that any domino tiling of
such a region can be obtained from any other tiling
by a sequence of moves,
each of which consists of applying a 90-degree rotation
to a pair of dominos that form a $2 \times 2$ block \cite{T}.
Such a move trades in an $a$-domino and a $b$-domino
for a $c$-domino and a $d$-domino, or vice versa.
The condition $ab=cd$ then guarantees that
the two tilings have the same weight.
Since such moves suffice to turn any tiling into any other,
all the tilings have the same weight,
and the probability distribution is uniform.
This property is called ``conditional uniformity,'' because if one
conditions on the tiling outside a simply-connected region, the
conditional distribution on tilings of the interior is uniform.

We do not propose any concrete interpretation for the weights $a,b,c,d$,
and it might be unreasonable to ask for one, given that the probability
distribution on matchings determined by $a,b,c,d$ is unaffected if all
four are multiplied by a constant.  There is a choice of scaling that
makes most of our formulas relatively simple, namely, the scaling that
makes the arcsines of the four weights add up to $\pi$, or equivalently,
the scaling that makes the product
$$
(a+b+c-d)(a+b-c+d)(a-b+c+d)(-a+b+c+d)
$$
equal to the product
$$
4(ab+cd)(ac+bd)(ad+bc)
$$
(this is most easily seen from Theorem~\ref{pa}).
However, it is worth mentioning at least one instance in which the
scaling $\sin^{-1} a + \sin^{-1} b + \sin^{-1} c + \sin^{-1} d = \pi$
does not give the simplest possible formula for $a,b,c,d$.
Specifically, consider the normalized Aztec diamond
with vertices at $(\pm 1, 0)$ and $(0, \pm 1)$.
The arctangent law of \cite{CEP} gives one way of expressing
how $p_a,p_b,p_c,p_d$ vary throughout the normalized diamond
(corresponding to the fact that the respective densities of
north-, south-, east-, and west-going dominos vary
throughout large Aztec diamonds).
However, an even more compact way of stating this dependence
is via the formulas
\begin{eqnarray*}
a & = & \sqrt{(1+y)^2 - x^2}, \\
b & = & \sqrt{(1-y)^2 - x^2}, \\
c & = & \sqrt{(1+x)^2 - y^2}, \\
d & = & \sqrt{(1-x)^2 - y^2}.
\end{eqnarray*}

\subsection{Sketch of proof and preliminary statement of results}
\label{sketchsec}

The strategy behind our proof of the main theorem is roughly as
follows.  In the first few (and more qualitative) sections of the
article (Sections~\ref{beginpartI} through \ref{endpartI}), we cover
the set of all domino tilings of the region by subsets which are balls
in the uniform metric on height functions (that is, each subset
consists of tilings whose height functions are approximately equal);
each ball is associated with an asymptotic height function $h$ from a
bounded subset of $\R^2$ into $\R$.  Appealing to quantitative results
proved later in more technical sections of the article
(Sections~\ref{partitionsection} through \ref{paproof}), as well as to
combinatorial arguments, we show that the logarithm of the cardinality
of a ball is approximated by the double integral in \eqref{supoverh}
times half the area of the region.
However, we also show that there is a unique ball
in our cover
for which the corresponding $h$ maximizes the double integral,
and that the contribution that this ball
makes to the total number of tilings swamps all the other contributions.

In Section~\ref{detsection}, we compute the ``partition function'' for
tilings
(that is, the sum of the weights of all the tilings).
This computation relies on Kasteleyn's original work \cite{Kast1}.
The next few sections of the article, on which the first few sections
depend, are close in spirit to Kasteleyn's original work on the dimer
model on the torus.  (Indeed, the local entropy
$\ent(s,t)$ is equal to the asymptotic
entropy for dimer covers of the $n\times n$ torus as $n$ goes to
infinity,
where the edges have been assigned new weights as
described above, favoring those tilings that have tilt near $(s,t)$.)
Using Kasteleyn's method of Pfaffians, we show that
the dimer model in the presence of weights $a,b,c,d$
exhibits a phase transition when any of the
four weights equals the sum of the other three, and our method of
analysis also requires that special attention be given to the case in
which $a=b$ and $c=d$ (which was analyzed by Kasteleyn).
The calculations are lengthy, but the end
results, obtained in the final sections of the paper, are satisfyingly
simple.
For direct formulas for $p_a,p_b,p_c,p_d$ in terms of
$s={\partial f}/{\partial x}$ and $t={\partial f}/{\partial y}$,
see (\ref{probdef}) below.

Here is a loose statement of our result.
See Theorems~\ref{main} and \ref{nbhdentropy}
for more precisely quantified statements,
and the remainder of Section~\ref{introsec}
for the relevant definitions.

\begin{thm}
\label{maintheorem}
Let $R^*$ be a region in $\R^2$ bounded by a piecewise smooth, simple
closed curve $\partial R^*$.  Let $h_b : \partial R^* \rightarrow \R$ be a
function which can
be extended to a function on $R^*$ with Lipschitz constant at most $2$
in the sup norm.  Let $f : R^* \rightarrow \R$ be the unique such
Lipschitz function
maximizing the entropy functional $\Ent(f)$ \textup{(}see equations
\eqref{entdef}
and \eqref{entfuncdef}\textup{)},
subject to $f |_{\partial R^*} = h_b$.
(See Section~\ref{beginpartI}~for the proof of the asserted uniqueness.)

Let $R$ be a lattice region that approximates $R^*$ when rescaled by a
factor of $1/n$, and whose normalized boundary height function
approximates $h_b$.  Then the normalized height function of a random
tiling of $R$ approximates $f$, with probability tending to $1$ as $n
\rightarrow \infty$.

Furthermore, let $g : R^* \rightarrow \R$ be any function satisfying the
Lipschitz condition.  Then the logarithm of the number of tilings of $R$
whose normalized height functions are near $g$ is $n^2(\Ent(g)+o(1))$.
\end{thm}

In his thesis \cite{H}, H\"offe derives (in a less rigorous manner)
expressions equivalent to our \eqref{Zformula} and \eqref{pa}.  In
\cite{DMB}, Destainville, Mosseri, and Bailly set up a similar
framework to our variational principle, but they use a quadratic
approximation to the entropy, and they do not supply rigorous proofs.
Readers of this article might
also wish to read
\cite{Kenyon2} and \cite{P2}.
These two articles cover some of the same ground as this one, but with
more of an emphasis on phenomena and less concern with proofs.

Before we can continue our discussion of the variational principle
and state the main result more precisely
(Subsection~\ref{vpsec}), we need to review some background on
entropy and height functions.

\subsection{Height functions}
\label{hfsec}

Height functions are a geometrical tool discovered
in individual cases by Bl\"ote and Hilhorst \cite{BH} and
Levitov \cite{Lev} and independently studied by Thurston \cite{T},
who situated the idea in a less ad hoc, more general context.
Given a (connected and) simply-connected region $R$ that can be tiled
by dominos, domino tilings of $R$ are in one-to-one correspondence
with height functions on the set of lattice points in $R$, defined up
to an additive constant.  The height function representation is useful
because the difference between the heights of two vertices encodes
non-local information about a tiling.  Here we will quickly summarize
the basic definitions and properties of height functions.  For a more
extensive discussion aimed at applications to random tilings, see
Subsections~6.1--6.3 and Section~8 of \cite{CEP}; for a
geometrical point of view, see \cite{T}.

Color the lattice squares in the plane alternately black and white,
like a checkerboard.  (It does not matter which of the two ways of
doing this is used.)  Define a horizontal domino to be {\it
north-going\/} or {\it south-going\/} according to whether its leftmost
square is white or black, and define a vertical domino to be {\it
west-going\/} or {\it east-going\/} according to whether its upper square
is white or black.  These names come from the domino shuffling
algorithm from \cite{EKLP}, which is the main combinatorial tool used
to study Aztec diamonds.  We will not use domino shuffling, but the
division of dominos into these four orientations will nevertheless be
important, as will the coloring in general.  (It is intended that
north-going dominos correspond to $a$-edges, south-going to $b$-edges,
east-going to $c$-edges, and west-going to $d$-edges.  The
geometrical terminology is more pleasant when we have no weights in
mind.)

Let $R$ be a lattice region (i.e., a connected region composed of
squares from the unit square lattice), and let $V$ be the set of
lattice points within $R$ or on its boundary.  We will always
assume that $R$ can be tiled by dominos in one or more ways.
A {\it height function\/} $h$ on $R$ is a function from $V$ to $\Z$
that satisfies the following two properties for
adjacent lattice points $u$ and $v$ on which it is defined.
(We consider two lattice points to be adjacent only if
the edge connecting them is contained in $R$.)  First, if
the edge from $u$ to $v$ has a black square on its left, then $h(v)$
equals $h(u)+1$ or $h(u)-3$.  Second, if the edge from $u$ to $v$ is
part of the boundary of $R$, then $|h(u)-h(v)|=1$.

Given a height function $h$ on $R$, consider the set of all pairs of
adjacent lattice points $u,v$ with $|h(u)-h(v)|=1$; then it is not
hard to check that the set of all dominos in $R$ that are bisected by
such an edge taken together constitute a tiling of $R$.

Given a domino tiling of a simply-connected region $R$, we can reverse
the process, and find a height function that leads to the tiling.
Such a height function always exists, but it is not quite unique,
because one can add
a constant (integer) value to it everywhere to get another such
height function.  To avoid this ambiguity, we fix the height of one
lattice point on the boundary of $R$.  Then height functions
satisfying this constraint are in one-to-one correspondence with
domino tilings.

{}From this point on, when we use the term {\it lattice region\/}, we will
always implicitly assume that one of the lattice points on the
boundary has a
specified
height, and when we talk about height functions
we will always assume that they satisfy this condition.  Notice that
when the region is simply-connected, our assumption takes on an
especially pleasant form, because it is then equivalent to fixing the
heights along the entire boundary.  (However, when the boundary is in
several pieces, the situation is more complicated.  In this paper, we
will typically assume that our regions are simply-connected, although
we will always mention that assumption when we make it.)

To define a height function for a domino tiling of an $n \times n$ torus,
view the torus as an $n \times n$ square with opposite sides identified,
and view this square as sitting centered
inside an $(n+2) \times (n+2)$ square.
A tiling of the torus determines a covering of the smaller square
by dominos that lie inside the larger square.  (Dominos that jump
from one side of the small square to the other correspond to two dominos
in the large square; the rest correspond to just one.)  Then define the
height function on the $n \times n$ square
as for a lattice region in the plane.  This height
function may not be well defined on the torus, since its values
can be different at two identified vertices, but this will not matter
for our purposes.  This definition is not really natural, but it is
convenient.  (See \cite{spaces} for a more natural definition of height
functions on tori.  We will not use their definitions or results.)

Let $R$ be a lattice region, and let $H(R)$ be the set of all height
functions on $R$.  We define a partial ordering $\le$ on $H(R)$ by
setting $h_1 \le h_2$ if $h_1(u) \le h_2(u)$ for every lattice point
$u \in R$.

The set $H(R)$ is not just a partially ordered set, but also a
lattice.  The join of two height functions $h_1$ and $h_2$ can be
defined by $(h_1 \hfsup h_2)(u) = \max(h_1(u),h_2(u))$, and their meet
by $(h_1 \hfinf h_2)(u) = \min(h_1(u),h_2(u))$.  To show that $H(R)$
is a lattice, we need to show that $h_1 \hfsup h_2$ and $h_1 \hfinf
h_2$ are height functions.  Notice that, at each lattice point $u$,
the value of height functions at $u$ is determined modulo $4$
(independently of the specific tiling); this assertion is trivial when
$u$ is the point at which we have fixed the values of height
functions, and if it is true at some particular lattice point, then
the definition of a height function immediately implies that it is
true at all adjacent lattice points.  Thus, if $h_1(u)$ is unequal to
$h_2(u)$, then they differ by at least $4$.  It now follows from
the definition of a height function that if $h_1(u) > h_2(u)$,
then $h_1(v) \ge h_2(v)$ for all lattice points $v$ adjacent
to $u$.  Hence, given any two adjacent lattice points,
$h_1 \hfsup h_2$ (or $h_1 \hfinf h_2$) agrees
with $h_1$ at both points, or agrees with $h_2$ at both of them,
which is what is needed to show that $h_1 \hfsup h_2$ and $h_1 \hfinf
h_2$ are height functions.  (See \cite{P1} for details and more
general results.)

As Figures~\ref{fig-square} and \ref{fig-aztec} demonstrate, the
precise boundary conditions of a region can have a dramatic effect on
the behavior of typical tilings.  Height functions provide the proper
tool for gauging the effect of boundary conditions on random tilings.
(Proposition~20 of \cite{CEP} is one way to make this claim precise.
It says, roughly, that the height function of a typical tiling depends
continuously on the heights on the boundary of the region.)
For example, one can compute the rate of change of the height
function in terms of the probabilities of finding dominos in
given locations.  If a region has statistically homogeneous random
tilings, then the typical height function should be nearly planar.
Because the boundary heights for Aztec diamonds are highly non-planar,
typical tilings cannot be homogeneous.

Define the {\it average height function\/} of a lattice region $R$ to be
the average of all height functions on $R$.  (Of course, it is almost
never a height function itself.)  Theorem~21 of \cite{CEP} implies
(roughly) that if $R$ is large, then almost all height functions on
$R$ approximate the average.  Thus, the problem of describing typical
tilings is reduced to the problem of describing the average height
function.

\subsection{Entropy}

The {\it entropy\/} of a random variable that takes on $n$ different
values, with probabilities $p_1,\dots,p_n$, is defined as $\sum_{i=1}^n
-p_i \log p_i$ (with $0 \log 0 = 0$).  For a uniformly distributed
random variable, the entropy is simply the logarithm of the number of
possible outcomes.  We will nevertheless
often
need to deal with entropy for a non-uniform distribution.
In general, we denote the entropy of a random variable $X$ by
$\ent(X)$;  there should be no confusion with our use of ``$\ent(\cdot)$'' to
denote local entropy depending on tilt.

We will sometimes speak of the {\it conditional entropy\/} of a
discrete random variable, conditional upon some event; this is
just the entropy of the conditional distribution determined by
the conditional probabilities.

The only fact about entropy that we will need other than the
definition is the following standard fact about
conditional entropy (the proof is straightforward).

\begin{lem}
\label{entropy}
Suppose $X$ is a random variable that takes on values
$x_1,\dots,x_n$, and suppose that $\{x_1,\dots,x_n\}$ is partitioned
into blocks $B_1,\dots,B_m$.  Let $B$ be the random variable that
tells which block $X$ is in, and let $X_i$ be the random variable that
takes on values in $B_i$ according to the conditional distribution of
$X$ given that $X \in B_i$.  Then the entropy of $X$  is given by
$$
\ent(X) = \ent(B) + \sum_{i=1}^m \Prob(x \in B_i) \ent(X_i).
$$
\end{lem}

When we deal with entropy for tilings of a lattice region $R$, we will
always normalize it by dividing by half the area of $R$, so that we
measure the information content per domino.  Thus, the normalized
entropy of a set of tilings of a lattice region $R$ (under the uniform
distribution) is the logarithm of the number of tilings in the set,
normalized by dividing by the number of dominos in a tiling of $R$.
When we refer to the entropy of a region $R$, we mean the entropy of
the set of all tilings of $R$, under the uniform distribution.

\subsection{The variational principle}
\label{vpsec}

Let $R^*$ be a region in $\R^2$ bounded by a piecewise smooth, simple
closed curve $\partial R^*$.  (We will always assume our curves do not
have cusps.)  Suppose $R$ is a large, simply-connected lattice region.
We can normalize by scaling the coordinates by a factor of $1/n$ (for
some appropriate choice of $n$). Suppose that the normalization of $R$
approximates $R^*$ (in a sense to be clarified below).  If we scale
the values of height functions by dividing their values by $n$, then
for large $n$ the {\it normalized height functions\/} that one obtains
approximate functions on $R^*$ that satisfy a Lipschitz condition with
constant $2$ for the sup norm.  The reason for this is that if $u$ and
$v$ are lattice points at sup norm distance at most $d$ from each
other within $R$ (i.e., they are connected by a path within $R$ of
length $d$, where the allowable steps are lattice edges or diagonal
steps) and $h$ is any height function, then one can check that
$|h(u)-h(v)| \le 2d+1$.  (One can prove this claim directly, or deduce
it from Lemmas~\ref{fournier} and \ref{supnorm}.)

Whenever we refer to a {\it 2-Lipschitz function\/} $f$ from a subset of
$\R^2$ to $\R$, we mean one that is locally Lipschitz with Lipschitz
constant $2$ for the sup norm, i.e., its domain is covered by open
balls in which every pair of points
$(x_1,y_1)$ and
$(x_2,y_2)$ satisfy
$$
|f(x_1,y_1)-f(x_2,y_2)| \le 2 \max(|x_1-x_2|,|y_1-y_2|).
$$
Notice that where such a function $f$ is differentiable, it must
satisfy
$$
|\partial f/\partial x| + |\partial f / \partial y| \le 2.
$$
Conversely, every differentiable function satisfying this condition is
$2$-Lipschitz.
We call $(\partial f/\partial x, \partial f/\partial y)$ the {\it
tilt\/} of the function (and use this terminology whether or not $f$ is
Lipschitz).  It is important to keep in mind Rademacher's theorem
(Theorem~3.16 of \cite{F}), which says that every Lipschitz function is
differentiable almost everywhere.

Suppose that $\partial R^*$ is a simple, closed, piecewise smooth
curve in $\R^2$ that bounds a region $R^*$.
We say that a function $h_b$ defined on
$\partial R^*$ is a {\sl boundary asymptotic height function\/}
if there exists a 2-Lipschitz function $h$ on $R^*$ such that
$h|_{\partial R^*} = h_b$, and we call such an $h$ an {\sl
asymptotic height function\/}.  Let $\Asymp(R^*,h_b)$ be the set
of all asymptotic height functions on $R^*$ with boundary
heights $h_b$;  notice that this set is convex and that it is
compact with respect to the sup norm.

Before continuing, we need to specify exactly what it means for one
region bounded by a simple closed curve to approximate another.
This
can be defined by
the Hausdorff metric on closed subsets of $\R^2$;  specifically,
two regions are defined to be within $\varepsilon$ of each other if
the $\varepsilon$-neighborhood of each one's boundary curve contains
the other's.

We also need to discuss what
approximation to within $\varepsilon$
means when there are boundary asymptotic height functions on the
curves.  Given regions $R_1$ and $R_2$ with boundary asymptotic height
functions $h_1$ and $h_2$, we say that $(R_1,h_1)$ is within
$\varepsilon$ of $(R_2,h_2)$ if for all $x_1 \in \partial R_1$ there
exists a $x_2 \in \partial R_2$ such that $d(x_1,x_2) < \varepsilon$ and
$|h_1(x_1) - h_2(x_2)| < \varepsilon$, and vice versa.

For technical reasons, it is most convenient to write out our
arguments in the case of approximation from within, where we have a
sequence $R_1,R_2,\dots$ of regions in the interior of $R^*$ whose
limit is $R^*$.  In most cases, this is easily arranged.  The regions
$R_1,R_2,\dots$ will be rescaled lattice regions, and if $R^*$ is a
star-shaped region, then by adjusting the scaling we can assume
approximation from within, without changing our asymptotic results.
However, when $R^*$ is not star-shaped, slightly more is required.
There are two ways to deal with such domains.  First, Proposition~20
of \cite{CEP} tells us that average height functions behave
continuously when boundary heights are perturbed, and this lets one
adjust regions so that their normalizations fit within the limiting
region, without changing the asymptotics.  Second, one can give a
direct proof by cutting a general region into star-shaped pieces.
Thus, we will be able to assume approximation from within in later
sections without loss of generality.

In the case of approximation from within, it will be convenient to use
a slightly different definition of $\varepsilon$-approximation.
If $R_2$ is in the interior of $R_1$, with boundary
height functions $h_1$ and $h_2$, then we say
that $(R_1,h_1)$ is within
$\varepsilon$ of $(R_2,h_2)$ if
their boundaries are within $\varepsilon$, and for every height
function $h$ on $R_1$ extending $h_1$, the restriction of $h$ to
$\partial R_2$ is within $\varepsilon$ of $h_2$.
This clearly generates the same topology as the definitions above, but
it is more convenient to work with, so we will use it in later sections.

We saw in the first paragraph of this subsection that if $R$
approximates $R^*$ when
rescaled, then normalized height functions on $R$ approximate
asymptotic height functions on $R^*$ (because of the
Lipschitz condition on
height functions).  Later, Proposition~\ref{hfexist} will tell us that
every asymptotic height function is nearly equal to a normalized
height function; that is, the class of asymptotic height functions
has not been defined too broadly.

We will show that the average height function on $R$ is determined by
finding the asymptotic height function $f$ that maximizes the integral
of local entropy, which depends on the tilt of $f$.

Given $(s,t)$ satisfying $|s|+|t| \le 2$ (i.e., a possible tilt for an
asymptotic height function), define the local entropy integrand
$\ent(s,t)$ as follows.  First define
%
%
\begin{eqnarray}
\label{probdef}
\notag p_a & = & \ \ \frac{t}{4}+\frac{1}{2\pi}\cos^{-1}\left(\frac{\cos(\pi
t/2)-\cos(\pi s/2)}{2}\right),\\
\label{Peqns}p_b & = & -\frac{t}{4}+\frac{1}{2\pi}\cos^{-1}\left(\frac{\cos(\pi
t/2)-\cos(\pi s/2)}{2}\right),\\
\notag p_c & = & -\frac{s}{4}+\frac{1}{2\pi}\cos^{-1}\left(\frac{\cos(\pi
s/2)-\cos(\pi t/2)}{2}\right), \quad \textup{and}\\
\notag p_d & = & \ \ \frac{s}{4}+\frac{1}{2\pi}\cos^{-1}\left(\frac{\cos(\pi
s/2)-\cos(\pi t/2)}{2}\right),
\end{eqnarray}
where the values of $\cos^{-1}$ are taken from $[0,\pi]$.
We
set $a=\sin(\pi p_a)$, $b=\sin(\pi p_b)$, $c=\sin(\pi p_c)$, and
$d = \sin(\pi p_d)$ when $|s|+|t|<2$,
and we define $\ent(s,t)$ as in \eqref{entdef}.

We will see later that there is good reason to believe that
the numbers $p_a,p_b,p_c,p_d$ correspond to the
local densities of the four orientations of dominos.  However, we do
not have a proof.

Define the entropy functional on $\Asymp(R^*,h_b)$ by
\begin{equation}
\label{entfuncdef}
\Ent(h) = \frac{1}{\area(R^*)}\iint_{R^*} \LocEnt{h} \, dx \,
dy.
\end{equation}
Suppose that $R_1,R_2,\dots$ is a sequence of simply-connected lattice
regions with specified boundary heights, such that when $R_n$ is
normalized by $n$ (or, say, a constant times $n$), as $n \rightarrow
\infty$ the boundary converges to
$\partial R^*$ and the boundary heights to $h_b$ (as specified above).
Without loss of generality, we can assume that the normalized
boundaries of $R_1,R_2,\dots$ all lie within $R^*$.
We know from Proposition~\ref{exists+unique} that there is a unique $f
\in \Asymp(R^*,h_b)$ that maximizes $\Ent(f)$.
We will prove (in Theorem~\ref{nbhdentropy}) that the entropy of
tilings of $R_n$ whose normalized height functions are close to $f$
(that is, the logarithm of the number of such tilings, divided by half
the area of $R_n$) is $\Ent(f)+o(1)$ as $n \to \infty$.

We can now prove a sharpened version of the claim made in
the second paragraph of Theorem~\ref{maintheorem}:

\begin{thm}
\label{main}
For each $\varepsilon > 0$, the probability that a normalized random
height function on $R_n$ differs anywhere by more than $\varepsilon$
from the entropy-maximizing function $f$ is exponentially small in
$n^2$ (i.e., is bounded above by $r^{n^2}$ for some $r<1$).
\end{thm}

\begin{proof}
Cover $\Asymp(R^*,h_b)$ with open sets around each asymptotic height
function $h$, such that the entropy for normalized height functions in
these sets is strictly less than $\Ent(f)$ unless $h=f$.
By
compactness, only finitely many of the sets are needed to cover
$\Asymp(R^*,h_b)$; we include the one corresponding to $h=f$.  Then
Theorem~\ref{nbhdentropy} implies that for $n$ large, the probability
that a random normalized height function on $R_n$ does not lie in the
open set around $f$ is exponentially small in $n^2$.
\end{proof}

Theorem~\ref{main} provides a much stronger large deviation estimate
than the best previous result (Theorem~21 of \cite{CEP}).  In
addition, it gives the first proof that under these conditions the
normalized average height function of $R_n$ converges to $f$ as $n
\rightarrow \infty$.  (It was previously not known to converge at all,
nor was a precise characterization of the entropy-maximizing function
known.)

It is worth pointing out that all our results generalize
to tilings with unit lozenges (as studied in, for example,
\cite{CLP}).  The combinatorial
results for lozenge tilings
(or, equivalently, the dimer model on a honeycomb graph)
are straightforward
modifications of those we present here for domino tilings;
the analytic results are special cases of ours (set one of
the four weights $a,b,c,d$
equal to $0$ to move from weighted domino tilings
of tori to weighted lozenge tilings).

In Sections~\ref{beginpartI} through~\ref{endpartI},
we will need the following three facts proved later
in the article:
\begin{enumerate}
\item For every tilt $(s,t)$ satisfying $|s|+|t| < 2$, there exist
unique weights $a,b,c,d$ (up to scaling) that satisfy $ab=cd$ and
give tilt $(s,t)$ (see Section~\ref{entasfnofslope}).

\item If an $n \times n$ torus has edge weights $a,b,c,d$ such
that $ab=cd$ and the tilt is $(s,t)$, then
the normalized entropy (for the probability distribution
on the tilings) is $\ent(s,t)+o(1)$ as $n \rightarrow
\infty$ (Theorem~\ref{ent}).

\item The local entropy function
$\ent(\cdot,\cdot)$ is strictly concave and continuous as a function of
tilt (Theorem~\ref{concavethm}).
\end{enumerate}

\section{Analytic Results on the Entropy Functional}
\label{beginpartI}

In this section, we will phrase all our results in terms of the local
entropy integrand, but they will depend only on its continuity and
strict concavity.  (Thus, if desired, the results can be applied in
different circumstances, for example to other sorts of tiling
problems.)  Everything in this section is fairly standard material
from geometric measure theory and the calculus of variations; for
example, Theorem~5.1.5 of \cite{F} is a more sophisticated version of
Lemma~\ref{semicontinuous}.  However, we will give complete proofs,
partly to make this part of the paper self-contained and accessible,
and partly because some of what we need does not seem to appear in the
literature in quite the form we would like.

In what follows, $h_b$ denotes a particular
boundary asymptotic height function
and $h$ ranges over the asymptotic height functions that restrict to $h_b$
on the boundary of the region $R$.

\begin{lem}
\label{semicontinuous}
The functional $\Ent : \Asymp(R^*,h_b) \rightarrow \R$, defined (as
above) by
$$
\Ent(h) = \frac{1}{\area(R^*)}\iint_{R^*}
\LocEnt{h} \,dx \,dy,
$$
is upper semicontinuous.
\end{lem}

For the proof of Lemma~\ref{semicontinuous}, as well as other
applications later in the paper, we will need to know how well we can
approximate an asymptotic height function $h$ by a piecewise linear
function.  The simplest way to achieve such an approximation is as
follows.  For $\ell > 0$, look at a mesh made up of equilateral
triangles of side length $\ell$ (which we call an $\ell$-mesh).
Consider any piecewise linear, 2-Lipschitz function $\tilde h$ that agrees
with $h$ at the vertices of the mesh and is linear on each triangle.
(Of course, $\tilde h$ depends on the mesh as well as on $h$.  It is
uniquely determined on each triangle that lies completely within
$R^*$, but not on those that extend over the boundary.)

\begin{lem}
\label{approx}
Let $h$ be an asymptotic height function, and let $\varepsilon > 0$.
If $\ell$ is sufficiently small then,
on at least a $1-\varepsilon$
fraction of the triangles in the $\ell$-mesh that intersect
$R^*$, we have the following two properties.  First, the
piecewise linear approximation $\tilde h$ agrees with $h$ to within
$\ell\varepsilon$.  Second, for at least a $1-\varepsilon$ fraction
(in measure) of the points $x$ of
the triangle, the tilt $h'(x)$ exists and is within $\varepsilon$ of
${\tilde h}'(x)$.
\end{lem}

Keep in mind that $h'(x)$ is the vector of partial derivatives of $h$
at $x$; it does not matter which norm we use to measure the distance
between tilts, but for the sake of specificity we choose the $\ell_2$
norm.
Notice that the second property implies that $\Ent(\tilde h) = \Ent(h)+o(1)$
(that is, $\Ent(\tilde h) \rightarrow \Ent(h)$ as $\varepsilon \rightarrow 0$).
(Keep in mind that $||h'||_1 \le 2$.)

\begin{proof}
We begin with the first of the two properties.  Recall that Lipschitz
functions are differentiable almost everywhere (Rademacher's theorem) \cite{F}.
For any point $x$ at which $h$ is differentiable, we have
$|h(x+d)-(h(x)+h'(x)\cdot d)| < \varepsilon |d|/2$ if $|d|$ is
sufficiently small, say $|d| \le r$ with $r > 0$ (where $r$ depends on
$x$).  If $x$ lies within an equilateral triangle of side length
$\ell$ with $\ell \le r$, then on that triangle we have the
approximation property we want, because there the two functions $d
\mapsto h(x+d)$ and $d \mapsto {\tilde h}(x+d)$ (the unique linear
function that agrees with $g$ on the corners of the triangle) both lie
within $\varepsilon \ell/2$ of $d \mapsto h(x)+h'(x)\cdot d$.

Given $\rho > 0$, let $S_{\rho}$ be the set of all $x$ such that $r$
(depending on $x$ as above) can be taken to be at least $\rho$.  Take
$\rho$ small enough that the measure of $S_{\rho}$ is at least
$(1-\varepsilon/3)$ times $\area(R^*)$.
(We can do that since $h$ is differentiable almost everywhere.)

Now take $\ell \le \rho$.  Look at any $\ell$-mesh, and the piecewise
linear approximation ${\tilde h}$ we get from it.  If $\ell$ is
sufficiently small, then all but a $o(1)$ fraction of the mesh
triangles lie entirely within the region.  At least a
$1-\varepsilon/3-o(1)$ fraction of them must intersect $S_{\rho}$,
which proves that the desired approximation property holds on at least
a $1-\varepsilon/3-o(1)$ fraction of the triangles.

For the second property, we will apply a result on metric density (see
Section~7.12 of \cite{Rudin}).  Let $U_1,\dots,U_n$ be open subsets
covering the set of possible tilts such that any two tilts contained
within the same subset differ by at most $\varepsilon$, and for $1 \le
i \le n$ let $V_i = \{x : h'(x) \in U_i \}$.  It follows from the
theorem on metric density that (for each $i$) if $\delta$ is
sufficiently small, then for all but an $\varepsilon/3$ fraction of
the points $x \in V_i$, at least a $1-\varepsilon$ fraction of the
ball of radius $\delta$ about $x$ lies in $V_i$ (and thus the tilt at
those points differs by at most $\varepsilon$ from $h'(x)$).  Now we
can take $\delta$ small enough that this result holds for all $i$ from
$1$ to $n$, and then as above it follows that for $\ell < \delta$, a
$1-\varepsilon/3-o(1)$ fraction of the mesh triangles in any
$\ell$-mesh lie entirely within $R^*$ and satisfy the second
property.

Thus, if $\ell$ is sufficiently small, at least a $1-\varepsilon$
fraction of the triangles satisfy both properties (since a
$1-\varepsilon/3-o(1)$ fraction satisfy each).
\end{proof}

\begin{lem}
\label{semielliptic}
Suppose that $\partial R^*$ is an equilateral triangle of side length
$\ell$, and
that the asymptotic height function $h$ satisfies $|h_b-\tilde h| <
\varepsilon \ell$ on $\partial R^*$, where $\tilde h$ is some linear
function.
Then
$$
\Ent(h) \le \Ent(\tilde h) + o(1)
$$
as $\varepsilon \rightarrow 0$.
\end{lem}

\begin{proof}
Because $\ent(\cdot)$ is concave,
\begin{equation}
\label{concaveineq}
\frac{1}{\area(R^*)}\iint_{R^*} \LocEnt{h} \, dx \, dy \le
\ent(h_{x,\av},h_{y,\av}),
\end{equation}
where
$$
h_{x,\av} = \frac{1}{\area({R^*})}\iint_{R^*} \frac{\partial
h}{\partial x} \,
dx \,dy
$$
and
$$
h_{y,\av} =
\frac{1}{\area({R^*})}\iint_{R^*} \frac{\partial h}{\partial
y}\,dx\,dy.
$$

We have
$$
(h_{x,\av},h_{y,\av}) = (\tilde h_{x,\av},\tilde h_{y,\av}) +
O(\varepsilon)
$$
(as one can see by computing the average by
integrating over cross sections and applying the fundamental theorem
of calculus), and hence
$$
\ent(h_{x,\av},h_{y,\av}) = \ent(\tilde h_{x,\av}, \tilde h_{y,\av}) +
o(1),
$$
since $\ent(s,t)$ is a continuous function of $(s,t)$.  Now combining
\eqref{concaveineq} with
$$
\Ent(\tilde h) = \area(R^*)\ent(\tilde h_{x,\av},\tilde h_{y,\av})
$$
yields the desired result.
\end{proof}

\begin{proof}[Proof of Lemma~\ref{semicontinuous}]
Let $h$ be any asymptotic height function, and consider the
neighborhood $U_\delta$ of $h$ consisting of all asymptotic height
functions within $\delta$ of $h$.  We need to show that given
$\varepsilon > 0$, if $\delta$ is sufficiently small, then for all $g
\in U_\delta$, $\Ent(g) \le \Ent(h)+\varepsilon$.

Let $\varepsilon' > 0$.  It follows from Lemma~\ref{approx} that if
$\ell$ is small enough, then the piecewise linear approximation
$\tilde h$ coming from an $\ell$-mesh satisfies $|h-\tilde h| <
\varepsilon' \ell$ on all but an $\varepsilon'$ fraction of the
triangles of the mesh, and $\Ent(\tilde h) = \Ent(h) + o(1)$ as
$\varepsilon' \rightarrow 0$.  Let $\delta < \varepsilon'\ell$.  Then
if $g \in U_\delta$, $|g-\tilde f| < 2\varepsilon' \ell$ on all but an
$\varepsilon'$ fraction of the triangles, and by
Lemma~\ref{semielliptic} the contribution to $\Ent(g)$ from the
non-exceptional triangles is at most $o(1)$ greater than the
corresponding contribution to $\Ent(\tilde h)$.  Of course, the
remaining $\varepsilon'$ fraction of the triangles contribute a total
of $O(\varepsilon')$.  Hence,
$$
\Ent(g) \le (1-O(\varepsilon'))\Ent(\tilde h)
 + o(1) +O(\varepsilon') = \Ent(h)+o(1).
$$
By choosing $\varepsilon'$ sufficiently small, one can make the $o(1)$
term less than $\varepsilon$.  Thus, $\Ent(\cdot)$ is upper semicontinuous.
\end{proof}

\begin{prop}
\label{exists+unique}
There is a unique asymptotic height function $f \in \Asymp(R^*,h_b)$
that maximizes $\Ent(f)$.
\end{prop}

\begin{proof}
For existence, we can use a compactness argument, since
$\Asymp(R^*,h_b)$ is compact.  Because the local entropy integrand is
bounded, $\Ent(\cdot)$ is bounded above, and we can choose a sequence
$h_1,h_2,\dots$ such that $\Ent(h_i)$ approaches the least upper bound
as $i \rightarrow \infty$.  By compactness, there is a subsequence
that converges, and by upper semicontinuity, the limit of the
subsequence must have maximal entropy.

Now uniqueness is easy.  Suppose that $f_1$ and $f_2$ are two
asymptotic height functions that maximize entropy, with $f_1 \ne f_2$.

Their derivatives cannot be equal almost everywhere, so
for some $\varepsilon > 0$ we must have
$$
\LocEnt{(f_1+f_2)/2} > \varepsilon +
\frac{\LocEnt{f_1}+\LocEnt{f_2}}{2}
$$
on a set of positive measure, by the strict concavity of $\ent(\cdot)$.  It
follows that
$$
\Ent\left(\frac{f_1+f_2}{2}\right) > \frac{\Ent(f_1)+\Ent(f_2)}{2},
$$
which contradicts the assumption that $\Ent(f_1)$ and $\Ent(f_2)$ were
maximal.  (Notice that since $f_1$ and $f_2$ are asymptotic height
functions, so is $(f_1+f_2)/2$.)  Therefore, only one asymptotic
height function can maximize entropy.
\end{proof}

\section{Combinatorial Results on Entropy}

In this section, we will compare the entropies of several regions with
nearly the same shape, but slightly different boundary conditions.  To
deal with this sort of situation, we will use partial height
functions.  Let $R$ be a simply-connected lattice region.  A {\it
partial height function\/} on $R$ is an integer-valued function $h$
defined on a subset of the lattice points in $R$, including all the
lattice points on the boundary, such that $h$ satisfies the following
condition.  If $u$ and $v$ are adjacent lattice points on which $h$ is
defined, such that the edge from $u$ to $v$ has a black square on its
left, then $h(v)-h(u)$ is $1$ or $-3$.  We call $h$ a {\it complete
height function\/} if it is defined on all the lattice points in $R$,
and we say that a complete height function is an {\it extension\/} of a
partial height function if they agree where both are defined.  We call
a partial height function a {\it boundary height function\/} if the set
of lattice points on which it is defined is the boundary of $R$ and
there exists a complete height function extending it.

Define a {\it free tiling\/} (or a ``tiling with free boundary conditions'')
of a region as a covering of the region
by dominos, no two of which overlap, and each of which contains at
least one cell belonging to the region.  A free tiling is just like
a tiling, except that the tiles are allowed to cross the boundary.
Given a boundary height function $h$, extensions of $h$ to $R$
correspond to tilings where dominos cross the boundary of $R$ at
certain specified places (namely the places where the height
changes by $\pm 3$).
In this way $h$ singles out a certain subset of the free tilings of $R$.
We could phrase all of our results in terms of free tilings,
but partial height functions will be a more convenient formulation.

Suppose we have a boundary height function on $R$, such that there is
an extension to a complete height function on $R$.  We can determine
the maximal extension $H$ (under the usual partial ordering) as
follows.  For any lattice point $v$ in $R$, look at boundary points
$w$ and paths $\pi$ joining $w$ to $v$ such that every edge of $\pi$
(oriented from $w$ towards $v$) has a black square on its left.  Call
such a path an {\it increasing path\/}, since it if does not cross any
dominos, then the height increases by $1$ after each edge.  Define a
{\it decreasing path\/} analogously.

\begin{lem}[Fournier]
\label{fournier}
For each lattice point $v \in R$, define $H_{\max}(v)$ as the minimum
of $h(w) + \length(\pi)$, where $w$ ranges over all boundary lattice
points and $\pi$ ranges over all increasing paths in $R$ from $w$ to
$v$, and define $H_{\min}(v)$ as the maximum of $h(w)-\length(\pi)$,
where this time $\pi$ ranges over all decreasing paths in $R$ from $w$
to $v$.  Then $H_{\max}$ is the maximal extension of $h$ to $R$, and
$H_{\min}$ is the minimal extension.
\end{lem}

For a proof, see \cite{Fournier}.  Notice that Lemma~\ref{fournier}
implies that if one changes the boundary height function $h$ by at
most a constant $c$ at each point, then $H_{\max}$ and $H_{\min}$
change by at most $c$ as well.

We can now prove a proposition we will need later, that connects
asymptotic height functions to actual height functions.  Suppose $R$
is a domino-tileable lattice region, such that rescaling $R$ by a
factor of $1/n$ gives a region whose boundary lies close to a region
$R^*$, where $\partial
R^*$ is a piecewise smooth, simple closed curve.  Without loss of
generality, we can suppose that the normalized region lies within
$R^*$ (as discussed in Subsection~\ref{vpsec}).

\begin{prop}
\label{hfexist}
Given an asymptotic height function $f$ which is within
$\varepsilon$ of the normalized boundary heights
of $R$, there is an actual height function on $R$ whose normalization
is within $\varepsilon +O(1/n)$ of $f$.
\end{prop}

\begin{proof}
Let $g$ be the largest height function on $R$ whose normalization
is less than or equal to $f$, i.e., the lattice sup of all height
functions below $f$.  (Technically, we must pick a lattice point and
restrict our attention to height functions that give it
height $0$ modulo $4$.  This ensures that all our
height functions are equal modulo $4$ and thus that the lattice
operations are well defined.)  Then all values of $g$ are within $8$
of what one would get simply by un-normalizing $f$, because if one
takes any lattice point, assigns it a value that is correct modulo $4$
and at least $4$ below the un-normalized value of $f$, and
looks at the minimal height function taking that value there,
then the Lipschitz constraint on $f$ implies that one stays below $f$.

The height function $g$ will typically not have the
correct boundary values on $R$.
Instead, it will have boundary
heights corresponding to some different boundary height function $b$,
although they will be within $\varepsilon n +O(1)$ of the actual boundary
heights for $R$.

To fix this problem, let $h$ be the minimal height function on $R$,
and let $H$ be the maximal one.  Then consider $(g \hfsup h) \hfinf
H$, which is a height function with the correct boundary values for
$R$.  It follows from Lemma~\ref{fournier} that $h$ and $H$ differ by
only $\varepsilon n+O(1)$ from the minimal and maximal extensions $h'$
and $H'$ of $b$, respectively.  Thus, since $(g \hfsup h') \hfinf H' =
g$, it follows that $(g \hfsup h) \hfinf H$ differs by only
$\varepsilon n +O(1)$ from $g$.  Because the normalization of $g$ differs
from $f$ by only $O(1/n)$, we see that
the normalization of $(g \hfsup h) \hfinf H$ differs from $f$ by
at most $\varepsilon +O(1/n)$, as desired.
\end{proof}

\begin{lem}
\label{supnorm}
The minimal increasing path length (with the path not restricted to
lie in any particular region) between two lattice points is always within
$1$ of twice the sup norm distance between them.
\end{lem}

\begin{proof}
Start at any lattice point and work outwards, labelling the other
points with the lengths of the shortest increasing paths to them.  It
is easy to prove by induction that on the square at sup norm distance
$n$ from the starting point, on two opposite sides the lengths
alternate between $2n$ and $2n+1$, and on the other two sides they
alternate between $2n$ and $2n-1$.
\end{proof}

We will prove the next three results in more generality than we need,
in order to make the precise hypotheses clear.  We will need to apply
them only to regions whose boundaries closely approximate a $k \times
k$ square or an equilateral triangle of side length $k$, such that the
heights along the boundary are nearly planar (in particular are fit by
a plane to within $\varepsilon k$ for some fixed $\varepsilon > 0$).

\begin{prop}
\label{planar}
Let $\varepsilon > 0$.  Suppose $R$ is a simply-connected lattice
region of diameter at most $n$ (i.e., every lattice point in $R$ is
connected to every other by a path within $R$ of length $n$ or less),
such that the heights on the boundary of $R$ are fit to within
$\varepsilon n$ by a plane with tilt $(s,t)$ satisfying $|s| + |t| \le 2$.
Then the average height function is given
to within $O(\varepsilon n)$ by that plane (if $n$ is sufficiently
large).
\end{prop}

Here ``$O(\varepsilon n)$'' simply denotes a quantity bounded in
absolute value by a fixed constant times $\varepsilon n$, for
sufficiently large $n$.

\begin{proof}
Consider a large torus, with edge weights $a,b,c,d$
chosen to give tilt $(s,t)$ and satisfying $ab=cd$ (see
Subsection~\ref{entasfnofslope}), and
view $R$ as being contained in the torus.
We will look at
random free tilings of $R$ generated according to the probability
distribution on weighted tilings of the torus.
If we fix the height of one point on the boundary of $R$,
then the average height function for these tilings
is given (exactly) by a linear function of the two position coordinates,
so that its graph is a plane.

If we select a random tiling
from this distribution,
then with probability
differing from $1$ by an exponentially small amount, the heights along
the border of the patch stay within $\varepsilon n$ of the plane, by
Proposition~22 of \cite{CEP}.  (It is not hard to check that all the
large deviation results from Section~6 of \cite{CEP}, such as
Propositions~20 and 22, apply to random tilings generated this way,
not just from the uniform measures on tilings of finite
regions.  All that matters is conditional uniformity, in the sense
that two tilings agreeing everywhere except on a subregion are
equally likely to occur; conditional uniformity
follows from $ab=cd$.)
Consider all possible boundary height functions on $R$
that stay within $\varepsilon n$ of a plane.
The average of the corresponding
average height functions, weighted by how likely they are to occur in
the weighted probability distribution,
is within $o(1)$ of the plane (for $n$ large).

By Proposition~20 of \cite{CEP}, all boundary height functions that
stay within $\varepsilon n$ of the plane have average height functions
within $O(\varepsilon n)$ of each other.  Since the average is within
$o(1)$ of the plane, each must be within $O(\varepsilon n)$ of it.
This completes the proof.
\end{proof}

\begin{lem}
\label{extremal}
Let $\varepsilon > 0$.  Suppose $R$ is a horizontally and vertically
convex lattice region of area $A$ with at most $n$ rows and columns,
such that $n \le \varepsilon A$.  Assume that the boundary heights
are fit to within $\varepsilon A/n$ by a plane with tilt $(s,t)$.
If $|s|+|t| \ge 2-\varepsilon$, then the entropy
of $R$ is $O(\varepsilon \log 1/\varepsilon)$, if $A$ is sufficiently
large.
\end{lem}

\begin{proof}
Without loss of generality, suppose that $s$ and $t$ are positive.

Consider any vertical line through $R$ on which the $x$-coordinate is
integral.  If the segment that lies within $R$ has length $k$, then
for every tiling of $R$, the difference between the number of
north-going dominos bisected by the line and the number of south-going
dominos bisected by it is $tk/4+O(\varepsilon A/n)+O(1)$, as one can
see by considering the total height change along the segment.
Similarly, on a horizontal segment of length $k$, the number of
west-going dominos bisected minus the number of east-going ones
bisected is $sk/4+O(\varepsilon A/n) + O(1)$.  If we add up all of
these quantities, the error term becomes $O(\varepsilon A) + O(n),$
which is $O(\varepsilon A)$ since $n \le \varepsilon A$.

The total number of dominos in any tiling of $R$ is $A/2$.  By the
results of the previous paragraph, this quantity is also twice the
number of east-going or south-going dominos plus $(s/4+t/4)A +
O(\varepsilon A)$.  We have $(s/4+t/4)A \ge ((2-\varepsilon)/4)A$.
Hence, the total number of east-going or south-going dominos is
$O(\varepsilon A)$.

Every tiling is determined by the locations of its east-going and
south-going dominos.  (To see why, recall that superimposing the
matchings corresponding to two tilings gives a collection of cycles.
Any disagreement between the two tilings yields a cycle of length at
least $4$, which must contain an east-going or south-going edge that
is in one tiling but not the other.)  Hence, there are at most
$$
O(\varepsilon A)\binom{A}{O(\varepsilon A)}2^{O(\varepsilon A)}
$$
tilings, since there are $O(\varepsilon A)$ possibilities for the
number of east-going or south-going dominos, at most
$$
\binom{A}{O(\varepsilon A)}
$$
places to put them, and at most
$$
2^{O(\varepsilon A)}
$$
ways to choose which are east-going.
It is not hard to check that Stirling's formula implies that
$$
\log \binom{A}{O(\varepsilon A)} =
O((\varepsilon\log 1/\varepsilon) A).
$$
Therefore, the entropy of $R$ is $O(\varepsilon\log
1/\varepsilon)$.
\end{proof}

\begin{prop}
\label{continuity}
Fix $k > 0$.  Let $\varepsilon > 0$, and let $R$ be a horizontally and
vertically convex region of area $A$ with at most $n$ rows and
columns, such that $kn^2 \le A$ and $n \le \varepsilon A$.  Suppose
$h$ is a boundary height function on $R$ that is fit to within
$\varepsilon A/n$ by a fixed plane.
Then for $A$ sufficiently large and $\varepsilon$
sufficiently small, the entropy of extensions of $h$ to $R$ is
independent of the precise boundary conditions, up to an error of
$O(\varepsilon^{1/2}\log 1/\varepsilon)$.
(The entropy does depend on the tilt of the plane, and this
proposition says nothing about whether it depends on the shape of $R$.
Notice also that the dependence on
$k$ is hidden within the implicit constant in the big-$O$ term.)
\end{prop}

\begin{proof}
Let $(s,t)$ be the tilt of the plane.  If $|s|+|t| \ge
2-\varepsilon^{1/2}$, then the conclusion follows from
Lemma~\ref{extremal}.  Thus, we can assume that the tilt satisfies
$|s|+|t| \le 2-\varepsilon^{1/2}$.

Suppose $g$ is another boundary height function on $R$, which agrees
with the same plane as $h$, to within $O(\varepsilon A/n)$.  We need
to show that extensions of $g$ and $h$ have nearly the same entropy.
Without loss of generality we can assume that $g \ge h$, since
otherwise we can go from $g$ to $g \hfsup h$ and from $h$ to $g \hfsup
h$.

Given any extension $f$ of $g$, let $H(f)$ be the infimum of $f$ and
the maximal extension of $h$, so that $H(f)$ is an extension of $h$.
Similarly, given any extension $f$ of $h$, let $G(f)$ be the supremum
of $f$ and the minimal extension of $g$, so that $G(f)$ is an
extension of $g$.

The maps $H$ and $G$ are not inverses of each other, but they come
fairly close to being inverses.  Given an extension $f$ of $h$,
$H(G(f))$ agrees with $f$ at every lattice point except those with
heights less than or equal to their heights in the minimal extension
of $g$.  By Lemma~\ref{fournier}, the minimal extension of $g$ is
within $O(\varepsilon A/n)$ of the minimal extension of $h$, so these
points have heights within $O(\varepsilon A/n)$ of their minimal
heights.  Similarly, given an extension $f$ of $g$, $G(H(f))$ agrees
with $f$ at all lattice points that are not within $O(\varepsilon
A/n)$ of their maximal heights.

By assumption, the tilt $(s,t)$ of the plane satisfies $|s|+|t| \le
2-\varepsilon^{1/2}$.  Thus, the height difference between any two
points in the plane is bounded by $2-\varepsilon^{1/2}$ times the sup
norm distance between them.  Therefore, Lemma~\ref{fournier} and
Lemma~\ref{supnorm} imply that the extreme heights at a point differ
from the heights on the plane by at least $\varepsilon^{1/2}$ times
the sup norm distance to the boundary.

Look at all lattice points not within sup norm distance
$(\varepsilon^{1/2} \log 1/\varepsilon) A/n$ of the boundary.  By
Proposition~\ref{planar}, at any such point the average height for
extensions of $g$ or $h$ is within $O(\varepsilon A/n)$ of the plane
(notice that because $kn^2 \le A \le n^2$, being $O(\varepsilon n)$ is
the same as being $O(\varepsilon A/n)$), and the extreme heights
differ from the plane by at least $(\varepsilon \log
1/\varepsilon)A/n$, so the extreme heights differ from the average
heights by an amount on the order of $(\varepsilon \log
1/\varepsilon)A/n$, for $\varepsilon$ small.
By Proposition~22 of \cite{CEP}, the probability that any such point
will have height within $O(\varepsilon A/n)$ of its extreme heights is
exponentially small in $n$.

Thus, given any point not within sup norm distance $(\varepsilon^{1/2}
\log 1/\varepsilon)A/n$ of the boundary, the probability that $H \circ
G$ or $G \circ H$ will not be the identity at that point is
exponentially small.  It follows that with probability nearly $1$, $H
\circ G$ and $G \circ H$ are the identity except on at most
$O((\varepsilon^{1/2} \log 1/\varepsilon)A)$ lattice points.  Thus,
the numbers of extensions of $g$ and $h$ differ by at most a factor of
$$
4^{O((\varepsilon^{1/2} \log 1/\varepsilon)A)},
$$
so the entropy of extensions of $g$ to $R$ differs from that of
extensions of $h$ by $O(\varepsilon^{1/2} \log 1/\varepsilon)$.
\end{proof}

\section{Proof of the Variational Principle}
\label{endpartI}

\begin{thm}
\label{precise}
Let $\varepsilon > 0$.  Suppose $R$ is an $n \times n$ square, with a
boundary height function $h$ fit to within $\varepsilon n$ by a plane
with tilt $(s,t)$ satisfying $|s|+|t| \le 2$.
Then for $n$ sufficiently large, the entropy of
extensions of $h$ to $R$ is
$$
\ent(s,t)+O(\varepsilon^{1/2} \log 1/\varepsilon),
$$
as is the entropy for free boundary conditions staying within
$\varepsilon n$ of the plane (i.e., the entropy for the set of all
free tilings of $R$ whose boundary heights stay within $\varepsilon n$
of the plane).
\end{thm}

\begin{proof}
We know from Proposition~\ref{continuity} that the entropy is
independent of the precise boundary conditions, but we still need to
prove that it equals $\ent(s,t)$.  To do so, we will compare with an
$n \times n$ torus that has edge weights $a,b,c,d$ satisfying $ab=cd$
and yielding tilt $(s,t)$.  (We can suppose that $|s|+|t|<2$, since
otherwise the result follows from Lemma~\ref{extremal} and
Proposition~\ref{continuity}.)
The torus is
obtained by identifying opposite sides of $R$, so that tilings of $R$
give tilings of the torus, but not vice versa.  Keep in mind that
because of the weighted edges in the torus, the probability
distribution on its tilings will not be uniform.  However, the
equation $ab=cd$ implies conditional uniformity (as mentioned
earlier), so if we fix the behavior on the boundary of $R$, then the
conditional distribution on extensions to the interior will be
uniform.

In Section~\ref{edgeprobsection}, we will define a set $W \subset
\N$ depending on $(s,t)$.  By Lemma~\ref{halfin}, for sufficiently
large even $n$, either $n$ or $n+2$ is in $W$.
In Proposition~\ref{ent1}, we show that the entropy of $n \times n$
tori with tilt $(s,t)$ converges to $\ent(s,t)$ as $n \rightarrow
\infty$ in $W$.

First we will suppose that $n \in W$, so that the entropy of the
$n \times n$ torus is $\ent(s,t)+o(1)$.

It follows from Proposition~22 of \cite{CEP} that with
probability exponentially close to $1$, in a random tiling of the
torus, the heights on the boundary of the square will be fit to within
$\varepsilon n$ by the average height function, which is a linear
function with tilt $(s,t)$.  The number of toroidal boundary
conditions is exponential in $n$, and by Proposition~\ref{continuity}
each has about the same entropy (except ones that are not nearly planar,
but they are
very unlikely to appear).  Lemma~\ref{entropy} tells us that the
entropy of the torus equals the average of the entropies for the
different boundary conditions, plus a negligible quantity for large
$n$ (since we have normalized by dividing by half of the area $n^2$).
Because
all the nearly planar boundary conditions have the same entropy, the
torus must as well, so since we know it has entropy $\ent(s,t)+o(1)$
as $n \rightarrow \infty$, each of the nearly planar boundary
conditions must have entropy $\ent(s,t) + O(\varepsilon^{1/2}
\log1/\varepsilon)$.  (Since $\varepsilon$ is fixed as $n \rightarrow
\infty$, we can absorb the $o(1)$ into the big-$O$ term.)

Now it is easy to deal with the case of $n \not\in W$.  If $n \not\in
W$, then $n+2 \in W$ and $n-2 \in W$.  The entropies for tilings
of $(n-2) \times (n-2)$, $n \times n$, and $(n+2) \times (n+2)$
squares with nearly
planar boundary conditions are nearly the same.  (To prove that,
embed an $(n-2) \times (n-2)$ square into an $n \times n$ one, and an
$n \times n$ one into an $(n+2) \times (n+2)$ one, extending the boundary
conditions arbitrarily.  Then the number of tilings increases with
each embedding, so the entropy of the $n \times n$ square is caught
between the other two, to within a $1+o(1)$ factor coming from the
differing areas.  By Proposition~\ref{continuity}, this result
holds for all nearly planar boundary conditions.)
Note that the same argument as above goes back from squares to the torus,
thus proving entropy convergence for all $n$ (not just $n \in W$).

The claim about free boundary conditions follows easily (since the
number of boundary conditions that stay within $\varepsilon n$ of the
plane is only exponential in $n$, and all of them have about the same
entropy).
\end{proof}

\begin{cor}
\label{triangle}
Theorem~\ref{precise} holds if $R$ is an equilateral triangle of
side length $n$, instead of a square.
\end{cor}

\begin{proof}
It is not hard to check that equilateral triangles satisfy
the hypotheses of Proposition~\ref{planar},
Lemma~\ref{extremal}, and Proposition~\ref{continuity}.
Thus, we just need to
deal with the case of a equilateral triangle whose boundary heights
are within $O(1)$ of being planar.

To prove that the entropy of the triangle is at least what we expect,
tile the triangle with smaller squares, such that their boundary
heights are within $O(1)$ of being planar.  (Of course, those near the
edges will stick out over the boundary, but if $n$ is large enough, we
can make the squares small enough compared to the triangle that only
an $\varepsilon$ fraction of the squares will cross the boundary.)
Except for an error of $O(\varepsilon)$ from the squares that cross
the boundary, the entropy of the triangle is at least the entropy of
the squares, which is what we wanted.

An analogous argument (involving tiling a square with triangles) shows
that the entropy of the triangle is at most what we expect.
\end{proof}

For the next theorem, we will use the same setup as in
Proposition~\ref{hfexist}.

\begin{thm}
\label{nbhdentropy}
Let $R^*$ be the region bounded by a simple closed curve, and let
$h_b$ be a boundary height function on $R^*$.
Suppose $R$ is a simply-connected lattice region
such that when $R$ is normalized by a factor of $1/n$, it
approximates $R^*$ to within $\delta$,
and its normalized boundary heights approximate the region $R^*$ with
boundary heights $h_b$ to within $\delta$;  we assume that the
normalization of $R$ lies within $R^*$.
Given an asymptotic
height function $h \in \Asymp(R^*,h_b)$, the
logarithm of the number of tilings of $R$
whose normalized height functions are within $O(\delta)$ of $h$ is
the area of $R$ times
$$
\Ent(h) + o(1)
$$
as $\delta \rightarrow 0$ (for $n$ sufficiently large).
\end{thm}

\begin{proof}
Notice that the set of tilings whose normalized height functions are
within $O(\delta)$ of $h$ is non-empty, by
Proposition~\ref{hfexist}. Call the set of such tilings $U_\delta$.

Fix $\varepsilon > 0$.  Choose $\ell$ small enough that we can apply
Lemma~\ref{approx} to the piecewise linear approximation $\tilde h$ to
$h$ derived from an $\ell$-mesh (with approximation tolerance
$\varepsilon$).  Then, as is pointed out after the
statement of Lemma~\ref{approx}, $\Ent(h) = \Ent(\tilde h) +
o(1)$.  We will take $\delta < \ell\varepsilon$, and show that the
entropy we want to compute is $\Ent(h)+o(1)$.

We know that $|h-\tilde h|<\ell\varepsilon$ on all but at most an
$\varepsilon$ fraction of the triangles in the mesh.  Those triangles
can change the entropy by only $O(\varepsilon)$, so we can ignore
them.  We can also ignore the $O(\delta)$ fraction of the triangles
that do not lie within the normalization of $R$ (which change the
entropy by $O(\delta)=O(\varepsilon)$).
We will call the triangles within the normalization of $R$ on which
$|h-\tilde h| < \ell\varepsilon$ the included triangles, and the
others the excluded triangles.

Let $g$ be any element of $U_\delta$.  The entropy of $U_\delta$ is
bounded below by the sum over all included triangles of the entropy of
$g$ restricted to that triangle, plus the $O(\varepsilon)$
contribution from the excluded triangles.  It is bounded above by the
same sum (including the $O(\varepsilon)$), but with free boundary
conditions on the included triangles (subject to the condition of
staying within $\ell\varepsilon$ of $f$).

We can now apply Corollary~\ref{triangle}.  It tells us that each
included triangle's contribution to the entropy of $U_\delta$ is
approximately equal to its contribution towards the entropy of $\tilde
h$.  It follows that our upper and lower bounds for the entropy of
$U_\delta$ both equal $\Ent(\tilde h) + O(\varepsilon) +
O(\varepsilon^{1/2}\log 1/\varepsilon)$. This gives us the desired
conclusion.
\end{proof}

Note that Theorem~\ref{nbhdentropy} implies Theorem~\ref{maintheorem}.

\newcommand{\PSbox}[3]{\mbox{\rule{0in}{#3}\includegraphics{#1}\hspace{#2}}}

\section{Overview of Remaining Sections}

In Section \ref{partitionsection} we compute the partition function
$Z_{n}(a,b,c,d)$ for matchings of the toroidal graph
$G_{n}=\Z^2/2n\Z^2$
with $4n^2$ vertices.  In Section
\ref{limitsection} we compute the limit, as $n\to\infty$, of
$Z_{n}^{1/(2n^2)}$.  In Section \ref{edgeprobsection} we compute the
limit of the edge-inclusion probabilities for edges of each type, with
respect to the measures $\mu_{n}$, and also a bound on the variance of
the number of edges of a fixed type in $G_{n}$.
This computation is only
done for $n$ in an infinite subset $W\subset\N$.  Since the variance
is $o(n^4)$, the measure is concentrating near tilings with the mean
number of edges of each type.  This fact allows us, in Section
\ref{entropysection}, to compute the limit for $n\in W$ of the
entropies.  As explained in the proof of Theorem~\ref{precise}, it
follows that the limit for arbitrary $n\to\infty$ must be the same as
the limit for $n\in W$.  In Section \ref{concavesection} we show that
the entropy is strictly concave as a function of the tilt $(s,t)$.
In Section \ref{PDEsection} we present the PDE which a $C^1$
entropy-maximizing Lipschitz function must satisfy
(at least in the distributional sense).

\section{The Partition Function}
\label{partitionsection}

Let the graph $G$ be the infinite square grid. Define an $a$-edge to
be a horizontal edge whose left vertex has even coordinate sum, a
$b$-edge to be a horizontal edge whose left vertex has odd coordinate
sum, a $c$-edge to be a vertical edge whose lower vertex has even
coordinate sum, and a $d$-edge to be a vertical edge whose lower
vertex has odd coordinate sum.  Let $a,b,c,d$ be four non-negative
real numbers.  Weight the $a$-edges with weight $a$, the $b$-edges
with weight $b$, and so on.  (For comparison with our earlier, more
geometrical terminology, $a$-edges are north-going, $b$-edges south-going,
$c$-edges east-going, and $d$-edges west-going;  in other words, points
with even coordinate sum correspond to white squares.)

We assume without loss of generality that $a\geq b$, $c\geq d$ and
$a\geq c$.

For $n$ an even positive integer let $G_{n}$
denote the quotient of $G$ by the action of translation
by $(2n,0)$ and $(0,2n)$. Then $G_{n}$ is a graph on the torus
and has $4n^2$ vertices and $8n^2$ edges ($2n^2$
edges of each type). A vertex $(a,b)$ of $G$ (where $a,b\in[0,2n-1]$)
will be denoted in what follows not by an ordered
pair $(a,b)$ but rather by an
ordered triple $(x,y,t)$, where $x=\lfloor a/2\rfloor$,
$y=\lfloor b/2\rfloor$, and
$t=1,2,3,$ or $4$ corresponding to $(a,b)$ being congruent to
$(0,0),(1,1),(1,0),$ or $(0,1)$ modulo $2$.

The partition function $Z_{n}$ is by definition the sum, over all perfect
matchings of $G_{n}$, of the product of the edge weights in the matching:
$$Z_{n}(a,b,c,d)=\sum_{\rm matchings} a^{N_a}b^{N_b}c^{N_c}d^{N_d},$$
where $N_a$ is the number of matched $a$-edges, etc.

There is a natural probability
measure $\mu_{n}=\mu_{n}(a,b,c,d)$  on the set of all matchings,
where the probability of a matching that has $N_a$ $a$-edges, etc.,
is
$$a^{N_a}b^{N_b}c^{N_c}d^{N_d}/Z_n.$$
The physical interpretation of $\mu$
is the following. Let $E_a$, $E_b$, $E_c$, and $E_d$
be energies associated with ``dimers'' on
an $a$-edge, $b$-edge, $c$-edge, and $d$-edge, respectively.
Define weights $a,b,c,d$
(also called ``activities'' in the statistical mechanics literature)
by
$$a=e^{-\beta E_a},~ b=e^{-\beta E_b},~ c=e^{-\beta E_c},~ d=e^{-\beta E_d},$$
where $\beta$ is a constant depending on the temperature.
Then $\mu$ is the
Boltzmann measure
associated to these energies; specifically,
the probability of a configuration of energy $E$
is proportional to $e^{-\beta E}$,
where $E$ is the sum of the energies of the individual dimers.

In what follows we will take $a,b,c,d$ as the fundamental quantities
and will not deal with $\beta$ or temperature as such.

Our concern will be with the situation in which $n$ goes to infinity
with the field-parameters $a,b,c,d$ fixed: the so-called ``thermodynamic
limit''.  Even though for each finite $n$, $Z_n(a,b,c,d)$ is a smooth
function (indeed a polynomial function) of the 4-tuple $(a,b,c,d)$,
the limit $Z(a,b,c,d) = \lim_{n \rightarrow \infty} Z_n^{1/(2n^2)}$
and other thermodynamic quantities need not be.
We will see that $Z(a,b,c,d)$ is $C^1$ everywhere but not $C^2$ in
the vicinity of the locus of $(a+b+c-d)(a+b-c+d)(a-b+c+d)(-a+b+c+d)=0$.

It's worth pointing out that the 4-parameter field determined by
$(a,b,c,d)$ actually only has two meaningful degrees of freedom:
one degree of freedom drops out because of the imposition of the
constraint $ab=cd$, and the other drops out by virtue of the fact
that multiplying $a,b,c,d$ by a constant has no effect on any of the
quantities of interest.
These two degrees of freedom correspond to the two degrees of freedom
associated with $(s,t)$ (the tilt).

\section{Determinants}
\label{detsection}
The goal of this section is to compute $Z_{n}$:
we use Proposition \ref{Kas} and equations
(\ref{detA})-(\ref{detA4}) below.

Given an enumeration of the $4n^2$ vertices of $G_n$,
we define the weighted adjacency matrix of $G_n$
as the $4n^2\times 4n^2$ matrix whose $i,j$th entry
is the weight of the edge connecting vertex $i$ to vertex $j$
(interpreted as 0 if there is no such edge).
Define a matrix $A_1$ by multiplying the weights on
vertical edges of the weighted
adjacency matrix of $G_{n}$ by $i=\sqrt{-1}$.
The matrix $A_2$ is obtained from $A_1$
by multiplying by $-1$ the weights on
the vertical edges from vertices $(j,n-1,4)$ to $(j,0,1)$
and edges from vertices $(j,n-1,2)$ to $(j,0,3)$ for all $j\in[0,n-1]$.
The matrix $A_3$ is obtained from $A_1$ by multiplying by $-1$ the weights on
horizontal edges from vertices $(n-1,k,2)$ to $(0,k,4)$
and horizontal edges from vertices $(n-1,k,3)$ to $(0,k,1)$,
for $k\in[0,n-1]$.
The matrix $A_4$ is obtained from $A_1$ by multiplying by $-1$ the weights on
both these sets of edges.
By the method of Kasteleyn \cite{Kast1}, we have the following proposition.
\begin{prop}\label{Kas}
For $i=1,2,3,4$ the quantities $\det A_i$ are
non-negative, satisfy $Z_{n}\geq \sqrt{\det{A_i}}$, and satisfy
$$Z_{n}(a,b,c,d
)=\frac12\left(-\sqrt{\det{A_1}}+\sqrt{\det{A_2}}+\sqrt{\det{A_3}}+
\sqrt{\det{A_4}}\right).$$
\end{prop}

\begin{figure}[htbp]
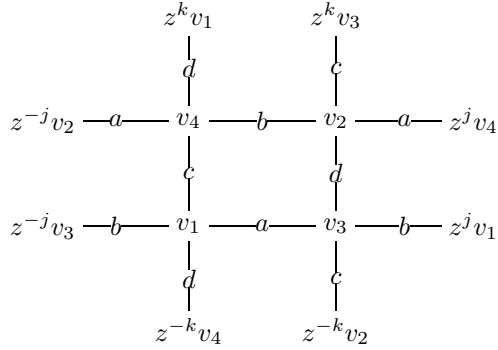

\begin{center}
\input evect.tex
\end{center}
\caption{\label{evects}A vector $v$ in the intersection of the
$z^j$-eigenspace of $T_{(2,0)}$ and the $z^k$-eigenspace of $T_{(0,2)}$.}
\end{figure}

Let $V$ denote the set of vertices of $G_{n}$.
The matrix $A_1$ operates on $\C^V$ in the following way:
for $f\in \C^V$ and $w\in V$,
$$(A_1f)_w = \sum_{v\in V}a_{vw}f_v.$$
Let $T_{(j,k)}$ be the linear operator on $\C^V$ corresponding to
the translation by $(j,k)$ on $G_{n}$.
The operators $T_{(2,0)}$ and $T_{(0,2)}$ commute with each other and
with $A_1$.  The eigenvalues of $T_{(2,0)}$
are $e^{2\pi ij/n}$ for integers $j\in[0,n-1]$.
The eigenspace of $T_{(2,0)}$
for eigenvalue $e^{2\pi ij/n}$ is $4n$-dimensional: a vector
$v$ in this eigenspace is determined by its coordinates in two consecutive
columns of $G_{n}$. Similarly $T_{(0,2)}$ has eigenvalues $e^{2\pi ik/n}$
and a vector in the $e^{2\pi ik/n}$-eigenspace is determined by its coordinates
in two consecutive rows of $G_{n}$.
The intersection of a maximal eigenspace of $T_{(2,0)}$ and one of $T_{(0,2)}$
is $4$-dimensional: a vector in the intersection is determined
by its coordinates on a $2\times 2$ square of vertices (as in Figure
\ref{evects}, vertices $v_1,v_2,v_3,v_4$).

Let $z=e^{2\pi i/n}$.  For $(j,k)\in[0,n-1]^2$
and $s\in\{1,2,3,4\}$ define a vector $e_{j,k}^{(s)}$ by
\begin{equation}
e_{j,k}^{(s)}(x,y,t)=
\begin{cases}
z^{jx+ky} & \textup{if $t=s$, and}\\
0 & \textup{otherwise}.
\end{cases}
\end{equation}
The $e^{(s)}_{j,k}$ are in the intersection of the $z^j$-eigenspace of $T_{(2,0)}$
and the $z^k$-eigenspace of $T_{(0,2)}$.
Let $S$ be the $4n^2\times 4n^2$ matrix whose columns are  these eigenvectors:
\begin{equation}
\label{S}
S=(e_{0,0}^{(1)},\dots,e_{0,0}^{(4)},
e_{1,0}^{(1)},\dots,e_{1,0}^{(4)},\dots,e_{n-1,n-1}^{(4)}).
\end{equation}
Note that $S$ satisfies $S^{-1}=\frac{1}{n^2}\overline{S^t}$,
where $^t$ denotes the transpose.

Because both $T_{(2,0)}$ and $T_{(0,2)}$ commute with $A_1$,
$S^{-1}A_1S$ has the block-diagonal form
\begin{equation}\label{SAS}S^{-1}A_1S=\left(\begin{array}{ccccc}
B_{0,0}&0&&&\\
0&B_{1,0}&0&&\\
&0&\ddots&0&\\
&&0&\ddots&0\\
&&&0&B_{n-1,n-1}\end{array}\right),
\end{equation} with $4\times 4$ blocks $B_{j,k}$ for
$j,k\in[0,n-1]$.
The block $B_{j,k}$ is the action of $A_1$ on the intersection of the
$z^j$-eigenspace of $T_{(2,0)}$ and the $z^k$-eigenspace of $T_{(0,2)}$.
We have
\begin{equation}\label{Bjk}
B_{j,k}=\left(\begin{array}{cccc}
0&0&a+bz^{-j}&i(c+dz^{-k})\\
0&0&i(d+cz^k)&b+az^j\\
a+bz^j&i(d+cz^{-k})&0&0\\
i(c+dz^k)&b+az^{-j}&0&0\end{array}\right)
\end{equation}
for the ordering $\{e_{j,k}^{(1)},\dots,e_{j,k}^{(4)}\}$
(see Figure \ref{evects}).

Since the upper right and lower left
$2$-by-$2$ subdeterminants of $B_{j,k}$ are
complex conjugates,
the determinant of $A_1$ is
$$\det A_1=\prod_{j=0}^{n-1}\prod_{k=0}^{n-1}
\left|2ab+a^2z^{-j}+b^2z^j+2cd+c^2z^{-k}+d^2z^k\right|^2,$$
so
\begin{equation}\label{detA}
\det A_1=\prod_{j=0}^{n-1}\prod_{k=0}^{n-1}
\left|\frac{(a+bz^j)^2}{z^j}+\frac{(c+dz^k)^2}{z^k}\right|^2.
\end{equation}

The matrices $A_2,A_3,A_4$ have similar determinants. For example
for $A_2$, let $R_{(0,2)}$ act on $\C^V$ by translation by $(0,2)$
followed by negation of the coordinates in the first two rows,
that is, vertices $(j,0,s)$ for $0\leq j\leq n-1$ and $s\in\{1,2,3,4\}$.
Then $R_{(0,2)}^n = -I$, and $R_{(0,2)}$ and $T_{(2,0)}$ commute
with $A_2$.
One thus finds that
\begin{equation}\label{detA2}
\det(A_2) =\prod_{j=0}^{n-1}\prod_{k=0}^{n-1}
\left|\frac{(a+be^{2\pi ij/n})^2}{e^{2\pi ij/n}}+\frac{(c+de^{\pi i(2k+1)/n})^2
}{e^{\pi i(2k+1)/n}}\right|^2,
\end{equation}
and similarly
\begin{eqnarray}\label{detA3}
\det(A_3)&=&\prod_{j=0}^{n-1}\prod_{k=0}^{n-1}
\left|\frac{(a+be^{\pi i(2j+1)/n})^2}{e^{\pi i(2j+1)/n}}+
\frac{(c+de^{2\pi ik/n})^2}{e^{2\pi ik/n}}\right|^2,\\
\label{detA4}
\det(A_4)&=&\prod_{j=0}^{n-1}\prod_{k=0}^{n-1}
\left|\frac{(a+be^{\pi i(2j+1)/n})^2}{e^{\pi i(2j+1)/n}}+
\frac{(c+de^{\pi i(2k+1)/n})^2}{e^{\pi i(2k+1)/n}}\right|^2.
\end{eqnarray}
We show that, with the exception of $A_1$, the square roots
of $\det A_i$ are polynomials in $a,b,c,d$.
Define
\begin{eqnarray}\label{P2}
P_2&=&\prod_{j=0}^{n-1}\prod_{k=0}^{n-1}
\left(\frac{(a+be^{2\pi ij/n})^2}{e^{2\pi ij/n}}+\frac{(c+de^{\pi i(2k+1)/n})^2
}{e^{\pi i(2k+1)/n}}\right),\\ \label{P3}
P_3&=&\prod_{j=0}^{n-1}\prod_{k=0}^{n-1}
\left(\frac{(a+be^{\pi i(2j+1)/n})^2}{e^{\pi i(2j+1)/n}}+\frac{(c+de^{2\pi ik/n})^2
}{e^{2\pi ik/n}}\right),\\ \label{P4}
P_4&=&\prod_{j=0}^{n-1}\prod_{k=0}^{n-1}
\left(\frac{(a+be^{\pi i(2j+1)/n})^2}{
e^{\pi i(2j+1)/n}}+\frac{(c+de^{\pi i(2k+1)/n})^2
}{e^{\pi i(2k+1)/n}}\right).
\end{eqnarray}
The functions $P_2,P_3$ and $P_4$ are polynomials in $a,b,c,d$.
Note that in (\ref{P2}), the involution $(j,k)\mapsto (-j,-k+1)$
maps each term to its complex conjugate (and this
involution has no fixed point).
Therefore $P_2=\sqrt{\det A_2}.$ Similarly $P_3=\sqrt{\det A_3}$
and $P_4=\sqrt{\det A_4}$.

Define
\begin{equation}\label{P1def}P_1=\pm\prod_{j=0}^{n-1}\prod_{k=0}^{n-1}
\left(\frac{(a+be^{2\pi ij/n})^2}{e^{2\pi ij/n}}+\frac{(c+de^{2\pi ik/n})^2
}{e^{2\pi ik/n}}\right),
\end{equation}
where the $+$ sign holds when $a<b+c+d$ and the $-$ sign holds
when $a>b+c+d$.
The involution $(j,k)\mapsto(-j,-k)$ maps each term in (\ref{P1def}) to
its complex conjugate; the product of the terms which are fixed under this
involution is
$$((a+b)^2+(c+d)^2)(-(a-b)^2+(c+d)^2)((a+b)^2-(c-d)^2)((a-b)^2+(c-d)^2).$$
This product
is positive or negative depending on whether $a<b+c+d$ or
$a>b+c+d$.  Thus on the domain $a<b+c+d$, $P_1$ is a polynomial
(taking non-negative values), and on the domain $a>b+c+d$,
$P_1$ is the negative of this polynomial (and this polynomial
also takes non-negative values).
Note that $P_1=\sqrt{\det A_1}.$ (The quantity $P_1$ is defined similarly
for all positive values of $a,b,c,d$; when one of $a,b,c,d$ is greater
than the sum of the other three, the $-$ sign is used, otherwise
the $+$ sign is used.)

{}From Proposition \ref{Kas} we have
\begin{equation}\label{Ps}
Z_{n}=\frac12(-P_1+P_2+P_3+P_4).
\end{equation}

\subsection{Eigenvalues and roots}
\label{roots}
Here we study in greater detail
the function $$q(z,w)=\frac{(a+bz)^2}z+\frac{(c+dw)^2}w,$$
where $a,b,c,d$ are non-negative reals.

Let
$$r(z)=cd+\frac{(a+bz)^2}{2z}.$$
Then
\begin{equation}
\label{q}
q(z,w)=\frac{c^2+2r(z)w + d^2 w^2}{w}=\frac{(dw-\alpha c)(dw-\beta c)}{w},
\end{equation}
where
\begin{equation}\label{alphabeta}
\alpha(z),\beta(z)=\frac{-r(z)}{cd}\pm\sqrt{\left(\frac{r(z)}{cd}\right)^2-1}.
\end{equation}
We will choose $\beta(z)$ to be the larger root (in modulus).

Let $\xi(z)=(a+bz)^2/z$ and $\eta(w)=-(c+dw)^2/w$, so that
$\xi(z)-\eta(w)=q(z,w)$. The critical points of $\xi(z),\eta(w)$ are
respectively $z=\pm a/b$ and $w=\pm c/d$.
Recall our assumption that $a\geq b$ and $c\geq d$.
\begin{lem}
\label{eta}
The map $\xi$ maps the punctured unit disk $\{z~:~|z|< 1,~z\neq 0\}$
injectively onto the exterior of the ellipse $E_1$
whose major axis has endpoints $-(a-b)^2,(a+b)^2$ and
whose minor axis has endpoints $2ab+i(a^2-b^2),2ab-i(a^2-b^2)$.
The ellipse $E_1$ has foci $0$ and $4ab$ and center
$2ab$.
In the case $a=b$, this ellipse degenerates to the line segment $[0,4ab]$.
The map $\eta$ maps the punctured unit disk $\{|w|< 1,~
w\neq 0\}$ injectively onto
the exterior of the ellipse $E_2$ with major axis $-(c+d)^2,(c-d)^2$,
minor axis $-2cd\pm i(c^2-d^2)$, and
foci at $0$ and $-4cd$.
\end{lem}
See Figures~\ref{ellipses} and~\ref{missellips}.

\begin{figure}[htbp]
\begin{center}
\input ell.pstex_t
\end{center}
\caption{\label{ellipses}The ellipses $E_1$ and
$E_2$ in the case $a<b+c+d$.}
\end{figure}

\begin{figure}[htbp]
\begin{center}
\input ell2.pstex_t
\end{center}
\caption{\label{missellips}The ellipses $E_1$ and $E_2$ in the case $a>b+c+d$.}
\end{figure}
\begin{proof}
Recall that $a\geq b$ by hypothesis. Suppose $a>b$.
We have $\xi(z)=a^2/z+b^2z+2ab$. For $|z|=1$ write $z=\cos t+i\sin{t}$.
Then
\begin{eqnarray*}\xi(z)&=&
a^2(\cos{t}-i\sin{t})+b^2(\cos{t}+i\sin{t})+2ab\\
&=&2ab+(a^2+b^2)\cos{t}- i(a^2-b^2)\sin{t}.
\end{eqnarray*}
Thus as $z$ runs counterclockwise around $S^1$,
$\xi(z)$ runs clockwise around the ellipse $E_1$.
Since $\xi$ has no critical points in the unit disk,
it is injective on the unit disk, mapping it to the exterior of $E_1$ (with $0$
mapping to $\infty$).

When $a=b$, $E_1$ degenerates to a segment; still,
$\xi$ maps the punctured
open unit disk injectively onto the exterior of the segment.

The case of $\eta$ is similar.
\end{proof}

The function $q$ can also be written
\begin{equation}\label{mod1}
q=\left(
\frac{a}{\sqrt{z}}+b\sqrt{z}+i\left(\frac{c}{\sqrt{w}}+d \sqrt{w}\right)\right)
\left(
\frac{a}{\sqrt{z}}+b\sqrt{z}-i\left(\frac{c}{\sqrt{w}}+d \sqrt{w}\right)\right).
\end{equation}
The coefficients of $a,b,c,d$ in either term of (\ref{mod1}) have modulus $1$ when
$|z|=1=|w|$.
Suppose $a>b+c+d$.
Then $q(e^{i\theta},e^{i\phi}) \neq 0$ for all $\theta$ and $\phi$, so
$E_1\cap E_2=\emptyset$, and hence $E_2$
is contained inside the bounded region delimited by $E_1$.
See Figure \ref{missellips}.
For each fixed $z$ satisfying $|z|=1$, $\xi(z)$ is on $E_1$ and
so by Lemma \ref{eta}
the roots $w$ of $0=q(z,w)=\xi(z)-\eta(w)$ are situated one outside and one
inside the unit disk.
Since $\beta$ is the larger of $\alpha$ and $\beta$,
$\alpha(z) c/d$ is the root inside the disk.

On the other hand if $a<b+c+d$ but not both $a=b$ and $c=d$,
then $(a-b)^2<(c+d)^2$ (and $(c-d)^2<(a+b)^2$
by the hypothesis that $a$ is the largest of $a,b,c,d$) so
the two ellipses intersect as in Figure~\ref{ellipses}
(because the places where they cross the $x$-axis
are interlaced). Let $(z_0,w_0)=(e^{i\theta_0},e^{i\phi_0})$ and
$(\overline{z_0},\overline{w_0})=(e^{-i\theta_0},e^{-i\phi_0})$
be the roots (satisfying $|z|=|w|=1$)
of $q(z,w)=0$, where $\theta_0\in(0,\pi)$
(the angle $\theta$ cannot be $0$ or $\pi$ since,
as we noted, the ellipses are interlaced on the $x$-axis).
Again by Lemma \ref{eta}, for each $z$ with $|z|=1$,
exactly one of the roots $w=\alpha c/d,\beta c/d$ is inside the
(closed) unit disk when
$-\theta_0\leq \theta\leq \theta_0$,
and for $\theta\not\in[-\theta_0,\theta_0]$,
both roots are outside.
We will again take $\alpha(z)c/d$ to be the smaller root.

In the case $a=b+c+d$, the two ellipses are tangent, and their
single intersection
point is at $z=-1,w=1$, that is, $(\theta_0,\phi_0)=(\pi,0)$.

In the case when both $a=b$ and $c=d$, the two degenerate ellipses
intersect only when $z=w=-1$, so that the single intersection
point is at $(\theta_0,\phi_0)=(\pi,\pi)$.

An important fact about the above four cases is that
$|\beta(z)c/d|>1$ always, and $|\alpha(z)c/d|>1$ unless
$\theta\in[-\theta_0,\theta_0]$.

\begin{figure}[htbp]
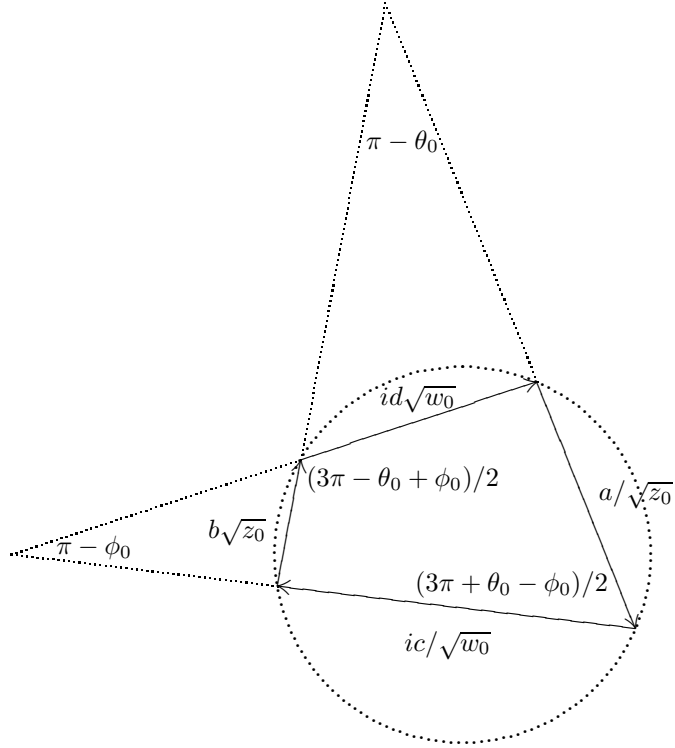

\begin{center}
\input cycquad.tex
\end{center}
\caption{\label{cycquad}The cyclic quadrilateral formed at
$(\theta_0,\phi_0)$.}
\end{figure}

Going back to the case $a\leq b+c+d$,
at the point $(\theta_0,\phi_0)$, the four quantities
$$\frac{a}{\sqrt{z_0}},
b\sqrt{z_0},\frac{ic}{\sqrt{w_0}},id\sqrt{w_0}$$ sum to zero by
(\ref{mod1}), assuming we choose the correct signs for $\sqrt{z_0}$ and $\sqrt{w_0}$.
They therefore form the edge vectors
of a quadrilateral. When taken in the order as in Figure \ref{cycquad},
the quadrilateral is in fact cyclic since
opposite angles sum to $\pi$:
The (interior) angle between sides $a/\sqrt{z_0}$ and $ic/\sqrt{w_0}$ is
$$\pi-\arg\left(\frac{a}{\sqrt{z_0}}\frac{\sqrt{w_0}}{ic}\right)
=\frac{3\pi}{2}-\arg
\sqrt{\frac{w_0}{z_0}},$$
and the (interior)
angle between sides $b\sqrt{z_0}$ and $id\sqrt{w_0}$ is
$$\pi-\arg\left(\frac{b\sqrt{z_0}}{id\sqrt{w_0}}\right)=
\frac{3\pi}{2}-\arg\sqrt{\frac{z_0}{w_0}}.$$
Summing these two gives angle $\pi$ (modulo $2\pi$, of course).

\subsection{The limit of the partition functions}
\label{limitsection}
\begin{thm}\label{Z} The limit
$$Z=Z(a,b,c,d)=\lim_{n\to\infty} Z_n^{1/(2n^2)}$$
exists, and
$$\log Z= \frac1{8\pi^2}
\int_0^{2\pi}\int_0^{2\pi}\log\left|\frac{(a+be^{i\theta})^2}{e^{i\theta}}+
\frac{(c+de^{i\phi})^2}{e^{i\phi}}\right| \,d\phi\, d\theta.$$
\end{thm}
Kasteleyn \cite{Kast1} computed $Z$ in the special case
$a=b$ and $c=d$.

\begin{proof}
{}From Proposition \ref{Kas}, for positive reals $a,b,c,d$ we have
$P_\ell\leq Z$ for all $\ell$. Combined with Proposition \ref{Kas}, this gives
$$\max_\ell P_\ell(a,b,c,d) \leq Z_{n}(a,b,c,d)\leq \frac32
\max_\ell P_\ell(a,b,c,d).$$
This implies that
$$\lim_{n\to\infty}Z_{n}^{1/(2n^2)}=\lim_{n\to\infty}
(\max_\ell P_\ell)^{1/(2n^2)}$$
(assuming these limits exist). Here note that
the $j$ for which $P_j=\max_\ell P_\ell$ may also depend on $n$.

Let $I$ denote the integral in the statement of the theorem.
We show that $({2n^2})^{-1}\log\max_\ell P_\ell$ converges to $I$.

Let
\begin{equation}\label{F}
F(\theta,\phi)=q(e^{i\theta},e^{i\phi})=
a^2e^{-i\theta}+2ab+b^2e^{i\theta}+c^2e^{-i\phi}+2cd+d^2e^{i\phi}.
\end{equation}
We then have $$P_1=\pm\prod_{j,k} F\left(\frac{2\pi j}n,
\frac{2\pi k}n\right),$$
and similar expressions for $P_2,P_3,$ and $P_4$
(see (\ref{P2})--(\ref{P1def})).
The quantity $n^{-2}\log(P_\ell)$ is a
Riemann sum for the integral of $\log|F|$.
The function
$\log|F|$ is continuous, hence Riemann integrable, on the complement of
a small neighborhood of its two (possible) singularities $(\theta_0,\phi_0)$
and $(-\theta_0,-\phi_0)$, defined in Section \ref{roots}.

We break the proof into four cases:
the case $a>b+c+d$, the case $a<b+c+d$ but not both $a=b$ and $c=d$,
the case $a=b$ and $c=d$,
and finally the case $a=b+c+d$.

If $a>b+c+d$,
then $\log|F|$ is continuous everywhere on $[0,2\pi]\times [0,2\pi]$
and so the Riemann sums converge to $I$.

In the case $a<b+c+d$ but not both $a=b$ and $c=d$,
there are two singularities.
It suffices to show that the Riemann sums $n^{-2}\log P_\ell$ for each
$\ell$
are small on a small
neighborhood of the singularities.  This is not quite true
since for some $\ell$ the product $P_\ell$ may have a factor in which
$(\theta,\phi)$ lands close to a singularity;
as a consequence this $P_\ell$ may be very small.
However, we will show that this can happen for at most one of the four
products $P_\ell$.

For each term in the four products (\ref{P2})--(\ref{P1def}),
$(\theta,\phi)$ is of the form $({j'\pi}/{n},{k'\pi}/{n})$,
for integers $j',k'$.  Furthermore for each pair of integers
$(j',k'),$ exactly
one of the four products has a term with $(\theta,\phi)=
({j'\pi}/{n},{k'\pi}/{n})$.
Thus at most one of the four products has a term closer than
$\pi/({2n})$ to the singularity $(\theta_0,\phi_0)$.
The {\em same} $P_\ell$ will have the term closest to the other singularity
$(-\theta_0,-\phi_0)$.

Fix a small constant $\delta>0$
and let $U_\delta\subset[0,2\pi]\times[0,2\pi]$
be the $\delta\times\delta$-neighborhood of a singularity.
In $U_\delta$ we use the Taylor expansion
\begin{equation}
\label{taylor}
F(\theta,\phi)= F_\theta(\theta_0,\phi_0)(\theta-\theta_0)+F_\phi(\theta_0,\phi_0)
(\phi-\phi_0)+ O((\theta-\theta_0)^2, (\phi-\phi_0)^2),
\end{equation}
and by Lemma \ref{notpar} below, the ratio of $F_\theta(\theta_0,\phi_0)$
and $F_\phi(\theta_0,\phi_0)$ is not real.

The sum of those
terms in $n^{-2}\log P_\ell$ for which $(\theta,\phi)\in U_\delta$
is
$$\frac1{n^2}\sum_{(\theta,\phi)\in U_\delta}\log|F(\theta,\phi)|=
\frac1{n^2}\left[
\sum_{(\theta,\phi)\in U_\delta}
\log\left|C_1(\theta-\theta_0)+C_2(\phi-\phi_0)\right|\right]
+ O(\delta^3),$$
for constants $C_1=F_\theta(\theta_0,\phi_0),C_2=F_\phi(\theta_0,\phi_0)$
with $C_1/C_2\not\in\R$. (We used here the fact that
$\log(x+O(x^2))=\log(x)+O(x)$ for small $x$;
we then could take the big-O term out of the
summation because of the $1/n^2$ factor.)

To bound this sum, note that
$|C_1x+C_2y|\geq C_3|x|+C_4|y|\geq C_5|x+iy|$ for some positive
constants $C_3,C_4,C_5$
since $C_1/C_2\not \in\R$.
We use polar coordinates around the singularity.
In the annulus around $(\theta_0,\phi_0)$ of inner radius $K/({2n})$ and
outer radius $(K+1)/(2n)$,
there are at most $\textup{constant}\cdot K$ points
$(\theta,\phi)$ which contribute to the sum,
and each such point contributes $\leq \left|n^{-2}
\log(\textup{constant}\cdot K/({2n}))\right|$ to the sum.
Therefore the sum  on $U_\delta$
(for those $P_j$ without terms within $\pi/(2n)$ of the singularity)
is bounded by
$$\frac{\textup{constant}
}{n^2}\sum_{1\leq K\leq \delta n}K\left|\log\frac Kn\right|=O(\delta^2\log\delta).$$

For the $P_j$ which does have a term closer than $\pi/(2n)$ to the singularity,
the above calculations give an {\em upper bound} on $({2n^2})^{-1}\log P_j$.
(Including in the factor close to the singularity only decreases the product.)
Thus we have shown that $({2n^2})^{-1}\log\max_\ell P_\ell$ converges to $I$.
This proves the convergence of $({2n^2})^{-1}\log Z$ to $I$.

In the case $a=b$ and $c=d$ we have only one singularity $(\pi,\pi)$,
and $P_1=0$. For the other $P_\ell$'s, when $(\theta,\phi)$ is close to
$(\pi,\pi)$ we have
\begin{eqnarray*}
F(\theta,\phi)&=&a^2(2+2\cos{\theta})+ c^2(2+2\cos{\phi})\\
&=&a^2(\theta-\pi)^2+c^2(\phi-\pi)^2 + O((\theta-\pi)^4,(\phi-\pi)^4).
\end{eqnarray*}
An argument similar to the previous case holds:
on $U_\delta$ we have
$$
\frac1{n^2} \sum_{(\theta,\phi)\in U_\delta} \log|F|=
\frac1{n^2}\left[
\sum_{(\theta,\phi)\in U_\delta}
\log\left|a^2(\theta-\pi)^2+c^2(\phi-\pi)^2\right|\right]
+ O(\delta^3).$$
Now $\log|a^2x^2+c^2y^2|\leq \log|C(x^2+y^2)|$.
Summing over annuli as before gives the bound.

Finally, suppose $a=b+c+d$. For any $\delta>0$ we have
$$Z_{n}(a+\delta,b,c,d)\geq Z_{n}(a,b,c,d)\geq Z_{n}(a-\delta,b,c,d)$$
(recall that coefficients of the polynomial $Z_{n}$ are non-negative).
For fixed $\delta$, the limits
$$\lim_{n\to\infty}Z_{n}(a\pm\delta,b,c,d)^{1/(2n^2)}$$
both exist, and (as we shall see in (\ref{Zformula}) below)
converge to the same value as $\delta\to 0$.
Thus
$$\lim_{n\to\infty}Z_{n}(a,b,c,d)^{1/(2n^2)}$$
exists and converges to this same value.
This completes the proof.
\end{proof}

The integral in Theorem~\ref{Z} can be written more usefully as follows.
As before
set $$r=r(z)=cd+\frac{(a+bz)^2}{2z}.$$
We can evaluate the first integral (the integral
with respect to $\phi$) in $I$ as follows.
We have (with $w=e^{i\phi}$)
$$\frac1{2\pi}\int_0^{2\pi}
\log\left|\frac{c^2+2rw+d^2w^2}{w}\right|d\phi =
\frac1{2\pi}\int_0^{2\pi}\log|(dw-\alpha c)(dw-\beta c)|\:d\phi$$
where $\alpha,\beta$ are chosen as in Subsection \ref{roots}.

Using the identity
$$
\frac1{2\pi}\int_0^{2\pi}\log|t+s e^{i\phi}|\:d\phi=
\begin{cases}
\log|t| & \textup{if $|t|>|s|$, and}\\
\log|s| & \textup{if $|s| \ge |t|$}
\end{cases}
$$
(note that the logarithmic singularity in the case $s=t$ makes
no contribution to the integral),
we find (with $z=e^{i\theta}$)
$$4\pi\log Z=
\int_{|c\alpha(z)|<d}\log d \:d\theta+\int_{|c\alpha(z)|>d}\log|c\alpha(z)|
\:d\theta$$
$$
+\int_{|c\beta(z)|<d}
\log d \:d\theta+\int_{|c\beta(z)|>d}\log|c\beta(z)|\:d\theta.$$

{}From Lemma \ref{eta} we know that $c|\beta(z)|>d$ for all $z$ on
the unit circle.
If $a\geq b+c+d$, then $|c\alpha(z)|\leq d$ for all $|z|=1$ and so
$$4\pi \log Z = 2\pi\log d+2\pi\log c
+\int_{-\pi}^{\pi}\log|\beta(z)|\:d\theta$$
or
\begin{equation}\label{a>bcd}
\log Z = \frac12\log d+\frac12\log c +\frac1{4\pi}
\int_{-\pi}^{\pi}\log|\beta(z)|\:d\theta.
\end{equation}

If $a<b+c+d$
then recall $|c\alpha(z)|<d$ if and only if
$\theta\in(-\theta_0,\theta_0)$.
Thus we have
$$4\pi\log Z =
\int_{-\theta_0}^{\theta_0}\log d \:d\theta+
\int_{\theta_0}^{2\pi-\theta_0}\log|c\alpha(z)|\:d\theta
+\int_{-\pi}^\pi\log|c\beta(z)|\:d\theta,$$
which gives (using $\alpha\beta=1$)
\begin{equation}\label{Zformula}
\log Z=\frac{\theta_0}{2\pi}\log d
+\left(1-\frac{\theta_0}{2\pi}\right)\log c +
\frac1{4\pi}\int_{-\theta_0}^{\theta_0}\log|\beta(z)| \:d\theta.
\end{equation}

Comparing (\ref{a>bcd}) and (\ref{Zformula}), we see that (\ref{Zformula})
holds for all $a,b,c,d$,
as long as we define $\theta_0=\pi$ when $a\geq b+c+d$.

\begin{lem}\label{notpar} For the function $F$ of (\ref{F}),
the ratio $F_\theta(\theta_0,\phi_0)/F_\phi(\theta_0,\phi_0)\not\in\R$
unless $a=b+c+d$ or both $a=b$ and $c=d$.
\end{lem}

\begin{proof}
We must show
\begin{equation}
\label{quotient}
\frac{-c^2e^{-i\phi_0}+d^2e^{i\phi_0}}{-a^2e^{-i\theta_0}+b^2e^{i\theta_0}}\not\in\R.
\end{equation}
For this proof only, let
$e^{i\theta}$ and $e^{i\phi}$ be square roots of $e^{i\theta_0},e^{i\phi_0}$,
respectively, with signs chosen so that
\begin{equation}\label{linear}
ae^{-i\theta}+be^{i\theta}+i(ce^{-i\phi}+de^{i\phi})=0
\end{equation}
(cf. (\ref{mod1})).
We can then factor (\ref{quotient}) as
$$\frac{(de^{i\phi}-ce^{-i\phi})}{be^{i\theta}-ae^{-i\theta}}\cdot
\frac{(de^{i\phi}+ce^{-i\phi})}{be^{i\theta}+ae^{-i\theta}},$$
and the second quotient is $i$.
That is, we are left to show that
$$
\textup{Im}\left(\frac{(d-c)i\cos{\phi}-(d+c)
\sin{\phi}}{(b-a)\cos{\theta}+i(b+a)\sin{\theta}}
\right)\neq0 .$$

Separating the real and imaginary parts of (\ref{linear}) we have
\begin{eqnarray}(a+b)\cos{\theta}+(c-d)\sin{\phi}&=&0, \textup{ and}\\
(-a+b)\sin{\theta}+(c+d)\cos{\phi}&=&0.
\end{eqnarray}
Solving these for $\sin{\phi},\cos{\phi}$ and plugging in to the above gives
$$
\textup{Im}\left(
\frac{(d-c)i\frac{a-b}{c+d}\sin{\theta}-(d+c)\frac{a+b}{d-c}\cos{\theta}}
{(b-a)\cos{\theta}+i(b+a)\sin{\theta}}\right)
$$ which is zero only if the real and imaginary parts of the numerator
and denominator are in proportion:
either $\sin\theta\cos\theta=0$ or (clearing denominators)
$$(d-c)^2(a-b)^2=(d+c)^2(a+b)^2.$$
Neither of these is possible (recall that $\theta_0=\pi$
only when $a=b$ and $c=d$), so the proof is complete.
\end{proof}

\section{The Edge-Inclusion Probabilities}
\label{edgeprobsection}
Let $\mu_{n}$ denote the measure on matchings of $G_{n}$,
where each matching has weight which is the product of the
edge weights of its matched edges.

The expected number of $a$-edges occurring
in a $\mu_{n}$-random matching is simply
$$\E(N_a)=\frac a{Z_{n}}
\frac{\partial Z_{n}}{\partial a}$$
(this follows from the definition of $\mu_{n}$).
{}From (\ref{Ps}), the probability $p_a(n)$ of a
particular $a$-edge
occurring in a $\mu_{n}$-randomly chosen matching is therefore
\begin{equation}\label{bigfrac}
p_a(n)=\frac{a}{2n^2}\frac{\frac{\partial}{\partial a}(-P_1+P_2+P_3+P_4)}
{-P_1+P_2+P_3+P_4}.
\end{equation}
{}From (\ref{P1def}) we obtain (assuming $P_1(a,b,c,d)\neq 0$)
\begin{equation}\label{ddaA}
\frac{\partial}{\partial a}P_1
=P_1\sum_{j,k}\frac{2(b+ae^{-2\pi ij/n})}{(a+be^{2\pi ij/n})^2e^{-2\pi ij/n}+
(c+de^{2\pi ik/n})^2e^{-2\pi ik/n}}.
\end{equation}
Note that this holds for all values of $a,b,c,d$
for which $P_1(a,b,c,d)\neq 0$, independently of whether or not $a>b+c+d$.
Similar expressions hold for $P_2,P_3,P_4$.

In what follows we may no longer assume
that $a$ is the largest of $b,c,d$
since we are computing a non-symmetric function $p_a$.
Recall that the quadruple $(a,b,c,d)$ determines
two possible singularities $\pm(\theta_0,\phi_0)$
of the function $F$ of (\ref{F}).

We will define a set $W\subset\N$, depending on $a,b,c,d$,
on which our remaining convergence arguments work.
If one of $a,b,c,d$ is greater than the sum of the others,
or if both $a=b$ and $c=d$, take $W=\N$.
If one of $a,b,c,d$ equals the sum of the others, we define $W$
below in the proof of Proposition~\ref{pamn}.
In the remaining case, there are two distinct singularities
$\pm(\theta_0,\phi_0)$, where $\theta_0\in(0,\pi)$.
If $\theta_0\neq\pi/2,$
let $W$ be the set of $n$ for which $\theta_0/\pi$
is not well approximated by rationals of denominator $n$,
in the following sense: for all integers $j$ we have
$$\left|\theta_0-\frac{\pi j}n\right|>\frac{1}{n^{3/2}}.$$
If $\theta_0=\pi/2$, define $W$ as above using $\phi_0$ instead:
note that $\theta_0$ and $\phi_0$ cannot both equal $\pi/2$,
for $F(\pi/2,\pi/2)=-i(a+bi)^2-i(c+di)^2$ cannot be zero
(its real part is $2ab+2cd$).

\begin{lem}
\label{halfin}
When $\theta_0\in(0,\pi)$,
for any sufficiently large even $n$, one of $n,n+2$ is in $W$.
\end{lem}

\begin{proof} Without loss of generality $\theta_0\neq\frac\pi2$.
Suppose
$$\left|\theta_0-\frac{\pi j}n\right|<\frac{1}{n^{3/2}}$$
and
$$\left|\theta_0-\frac{\pi j'}{n+2}\right|<\frac{1}{(n+2)^{3/2}}.$$
Then $j'$ must be equal to one of $j,j+1$ or $j+2$; but then
$$\left|\frac{\pi j}n-\frac{\pi j'}{n+2}\right|=
\frac{\pi\min\{2j,|2j-n|,|2j-2n|\}}{n(n+2)}\geq \frac{\textup{constant}}{n},$$
a contradiction.
\end{proof}
\begin{prop}
\label{pamn}
If none of $a,b,c,d$ equals the sum of the others, then
for $n$ tending to $\infty$ in $W$,
the edge-inclusion probability $p_a(n)$ converges to
\begin{equation}\label{sumintegral}p_a=
\frac{a}{4\pi^2}\int_0^{2\pi}\int_0^{2\pi}
\frac{(b+ae^{-i\theta})\:d\phi\:d\theta}{(a+be^{i\theta})^2
e^{-i\theta}+(c+de^{i\phi})^2e^{-i\phi}}.
\end{equation}
If one of $a,b,c,d$ equals the sum of the other three,
then $p_a(n)$ converges to the above integral along a subsequence $W$
containing at least one of each pair $n,n+2$ with $n$ even.
\end{prop}

\begin{proof}
We will show that the sum
\begin{equation}\label{psum}
\frac{a}{2n^2}\frac{\p}{\p a}\log P_1=\frac{a}{2n^2}\sum_{\theta,\phi}
\frac{2(b+ae^{-i\theta})}{F(\theta,\phi)}
\end{equation}
converges to the desired integral (\ref{sumintegral}), where the sum
is over $(\theta,\phi)=({2\pi j}/n,{2\pi k}/n)$.
Similar arguments hold for $P_2$, $P_3$, and $P_4$.
The value (\ref{bigfrac}) is then a weighted average of
these sums, with weights $\pm{P_\ell}/({2Z_{n}})$.
Since $Z_{n}\geq P_\ell\geq 0$ (see Proposition \ref{Kas}),
the weights are bounded in absolute value
(less than $1/2$) and sum to $1$.
Therefore the weighted average also converges to (\ref{sumintegral}).

We separate the proof into four cases.
In the first case,
where one of $a,b,c,d$ is greater than the sum of the other three,
there are no singularities ($F(\theta,\phi)$ is never zero),
so the summand is a continuous function on $[0,2\pi]^2$.
Therefore (\ref{psum}) converges to the integral in (\ref{sumintegral}).

For the second case,
suppose each of $a,b,c,d$ is strictly less than the sum of the others,
but we are not in the case where both $a=b$ and $c=d$.
Since $n\in W$, none of the four $P_\ell$
can have a term with $(\theta,\phi)$ within ${n^{-3/2}}$
of a singularity.  We claim that the sum (\ref{psum}) converges to the
integral (\ref{sumintegral}).
This is proved in the same manner as in Theorem~\ref{Z}: one needs only check
that the contribution on a small neighborhood of the singularities is small.
As before, let $U_\delta$ be a $\delta\times\delta$-neighborhood
of a singularity.
Ignore for a moment the single term closest to the singularity.
Lemma \ref{notpar} and (\ref{taylor})
give us the estimate
\begin{equation}\label{est}
\left|\frac{a}{2n^2}\sum_{(\theta,\phi)
\in U_\delta}\frac{2(b+ae^{-i\theta})}{F(\theta,\phi)}\right|
\end{equation}
$$
\leq\frac{2a}{2n^2}
\sum_{U_\delta}\frac{|b+ae^{-i\theta}|}{|C_1(\theta-\theta_0)+
C_2(\phi-\phi_0)+O((\theta-\theta_0)^2,
(\phi-\phi_0)^2)|},
$$
for constants $C_1=F_\theta(\theta_0,\phi_0),C_2=F_\phi(\theta_0,\phi_0)$
where $C_1/C_2\not\in\R$.
Now $$|C_1(\theta-\theta_0)+C_2(\phi-\phi_0)|\geq C_3|(\theta-\theta_0)
+i(\phi-\phi_0)|\geq C_4\max\{|\theta-\theta_0|,|\phi-\phi_0|\}$$
for some positive constants $C_3,C_4$.
Using
$1/({x+O(x^2)})=1/x+O(1)$ for small $x$, where $x=C_3|(\theta-\theta_0)
+i(\phi-\phi_0)|,$
we get
the bound
$$\frac{2a}{2n^2}
\sum_{U_\delta}\left[\frac{|b+ae^{-i\theta}|}{|C_1(\theta-\theta_0)+
C_2(\phi-\phi_0)|}+O(1)\right].$$
Taking the $O(1)$ term out of the summation turns it into a $O(\delta^2)$.

Summing over annuli
concentric about the singularity, we may
bound the left-hand side of (\ref{est}) by
$$O(\delta^2)+\frac{\textup{constant}}{n^2}\sum_{1\leq K\leq \delta n}K\cdot
\frac{n}{K}=O(\delta).$$
The single term closest to the singularity contributes a negligible
amount
$$\frac{\textup{constant}}{n^2}\cdot n^{3/2}.$$

In the third case, when both $a=b$ and $c=d$, we have
$(\theta_0,\phi_0)=(\pi,\pi)$.
Then $P_1=0$ and the other $P_\ell$ are nonzero. Furthermore the pairs
$(\theta,\phi)$ appearing in the products for $P_2,P_3,P_4$ do not come
within distance $\pi/n$ of the singularity.
Since near $(a,b,c,d)$, $P_1$ is a polynomial taking non-negative values
which is zero at $(a,b,c,d)$,
it must have a double root there. Thus its derivative with respect to $a$
is zero at $(a,b,c,d)$. We can therefore remove $P_1$ from the expression
(\ref{bigfrac}) and just deal with the remaining three $P_\ell$.
When $a=b$ and $c=d$, (\ref{est}) becomes
\begin{eqnarray*}
\frac{2a}{2n^2}\left|\sum_{U_\delta}
\frac{a(1+e^{-i\theta})}{a^2(2+2\cos\theta)+c^2(2+2\cos\phi)}\right|
&\le&
\frac{2a^2}{2n^2}\left[\sum_{U_\delta}
\frac{|i(\theta-\pi)|}{|a^2(\theta-\pi)^2+c^2(\phi-\pi)^2|}\right]\\
&& \phantom{}+ O(\delta^3),
\end{eqnarray*}
where we used $1/({x^2+O(x^4)})=1/{x^2}+O(1)$.
This is $O(\delta)$ (sum over annuli
as before).

Finally we consider the case when one of $a,b,c,d$ is equal to the sum of the
other three. Suppose first that $a=b+c+d$.
Note that $p_a(n)(a,b,c,d)$ is a monotonic increasing function of $a$:
the expected number of $a$-edges increases with their relative weight.
For each $\delta>0$ sufficiently small,  choose $n$ so that
$$
|p_a(n)(a-\delta,b,c,d)-p_a(a-\delta,b,c,d)|<\delta.
$$
Such an $n$ exists because $(a-\delta,b,c,d)$ is in the domain
of case two, above. By monotonicity
$$p_a(n)(a,b,c,d)\geq p_a(n)(a-\delta,b,c,d)\geq
p_a(a-\delta,b,c,d)-\delta$$ for this $n$.
Take a sequence of $\delta$'s tending to $0$. On the corresponding sequence
of $n$'s,
$p_a(n)(a,b,c,d)$ tends to $1$,
which is equal to the value of the integral $p_a(a,b,c,d)$ (see
Theorem~\ref{pa} below).
The set $W$ in this case is obtained from the concatenation of the appropriate
subintervals of the sets $W$
that are defined for each quadruple $(a-\delta,b,c,d)$.

Since $p_a(n)(a,b,c,d)\to 1$ we must have that
\begin{eqnarray*}p_b(n)(a,b,c,d)&\to& 0,\\
p_c(n)(a,b,c,d)&\to& 0,\\
p_d(n)(a,b,c,d)&\to& 0.
\end{eqnarray*}
This (and symmetry) takes care of the remaining cases.
\end{proof}

A rather lengthy calculation yields the following result:
\begin{thm}\label{pa} If $a\geq b+c+d$, then $p_a=1$.
If one of $b,c,d$ is larger than the sum of the other three of $\{a,b,c,d\}$,
then $p_a=0$. Otherwise,
let $Q$ be a cyclic quadrilateral with edge lengths $a,c,b,d$
in cyclic order.  Then $p_a$
is $1/(2\pi)$ times the angle of the arc cut off by the edge $a$ of $Q$.
That is,
\begin{equation}\label{paeqn}
p_a=\frac1\pi\textup{sin}^{-1}\left(\frac{a\sqrt{(a+b+c-d)(a+b-c+d)(a-b+c+d)(-a+b+c+d)}}{2
\sqrt{(ab+cd)(ac+bd)(ad+bc)}}\right).
\end{equation}
(Here the branch of the arcsine that is needed is the one given by
the preceding geometrical condition.)
\end{thm}

The proof is in Section \ref{paproof}.
Note that by this theorem, in the case of non-extremal tilt,
the $4$-tuple
$$(\sin(\pi p_a),\sin(\pi p_b),
\sin(\pi p_c),\sin(\pi p_d))$$ is proportional to $(a,b,c,d)$.
Since a constant of proportionality has no effect on the measures
$\mu_{n}$, we can assume that $a=\sin(\pi p_a)$, etc.

We also note a simple relation between the singularity
$(\theta_0,\phi_0)$
and the edge-inclusion probabilities
(hereafter simply called edge probabilities),
which will be useful later.
{}From Figure \ref{cycquad} and Theorem~\ref{pa},
when $a\geq b$ and $c\geq d$ we find
\begin{eqnarray}\label{thetapcpd}
\theta_0 &=& \pi-\pi(p_c-p_d), \textup{ and}\\
\phi_0 &=& \pi-\pi(p_a-p_b).
\end{eqnarray}

We now bound the variance of $N_a$, the number of $a$-edges in
a matching.

\begin{prop}
\label{sigmasquared}
For all $a,b,c,d$ and
for $n\in W$, $\sigma^2(N_a)=o(n^4)$.
\label{varbd}
\end{prop}

\begin{proof}
The graph $G_{n}$ has $2n^2$ $a$-edges. For $k\in[1,2n^2]$ let $q_k$
be the $\{0,1\}$-valued random variable indicating the presence of the $k$th
$a$-edge in a random matching.
Then $N_a=q_1+\dots+q_{2n^2}$, and so
\begin{equation}\label{var}
\sigma^2(N_a)=\sum_k \sigma^2(q_k)+\sum_{k\neq\ell}(\E(q_kq_\ell)-
\E(q_k)\E(q_\ell)).
\end{equation}

We have $\E(q_k)=p_a(n)$ and so $\sigma^2(q_k)=p_a(n)-p_a(n)^2$.
In the case when one of $a,b,c,d$ is greater than the sum
of the others, we know from Proposition \ref{pamn}
that $p_a$ converges to $1$ or $0$; as a consequence
$\sigma^2(q_k)\to 0$, and the covariances converge to $0$ also,
so $\sigma^2(N_a)=o(n^4)$ as well.
Similarly when one of $a,b,c,d$ equals the sum of the others; then
$p_a$ converges to $1$ or $0$ along $W$, and so $\sigma^2(N_a)$
is $o(n^4)$ on this same subsequence.

The remaining cases require more work.
By (a straightforward extension of) Theorem~6 of \cite{Kenyon1}, we have
\begin{equation}\label{detavg}
\E(q_kq_\ell)=
\frac{-|(A_1^{-1})_{q_k,q_\ell}|P_1+|(A_2^{-1})_{q_k,q_\ell}|P_2
+|(A_3^{-1})_{q_k,q_\ell}|P_3+|(A_4^{-1})_{q_k,q_\ell}|P_4}{-P_1+P_2+P_3+P_4},
\end{equation}
where
$$|(A_i^{-1})_{q_k,q_\ell}|=\det\left(\begin{array}{cc}A_i^{-1}(v_k,w_k)&
A_i^{-1}(v_k,w_\ell)\\A_i^{-1}(v_\ell,w_k)&A^{-1}(v_\ell,w_\ell)\end{array}\right),$$
and where $v_k,w_k$ are the vertices of the edge associated with $q_k$
($v_k$ being the left vertex) and $v_\ell,w_\ell$ the vertices of the
edge associated with $q_\ell$ ($v_\ell$ being the right vertex).  The
inverses $A_i^{-1}$ are only defined when the corresponding $P_i$ are nonzero.

When $n\in W$ tends to $\infty$ the diagonal entries
$A_i^{-1}(v_k,w_k),~A_i^{-1}(v_\ell,w_\ell)$
tend to $p_a$ (see Proposition \ref{pamn}). Writing each $2\times 2$ determinant
in (\ref{detavg}) as the product of the two diagonal entries minus the product
of the off-diagonal entries, the diagonal entries can be taken out of
the quotient in (\ref{detavg}) and contribute $p_a^2+o(1)$.
It remains to estimate the contribution of the off-diagonal entries.

We can compute the inverses of the $A_i$ as follows.
Note that from (\ref{Bjk}) we have
$$
B_{j,k}^{-1}=\left(\begin{array}{cc}0&D_1\\D_2&0\end{array}\right),$$
where
$$D_1=\frac{1}{(a+bz^j)^2z^{-j}+(c+dw^k)w^{-k}}\left(\begin{array}{cc}
b+az^{-j}&-i(d+cw^{-k})\\
-i(c+dw^k)&a+bz^j\end{array}\right)$$
(we won't need the expression for $D_2$).
Recall the definition of the matrix $S$ of (\ref{S}).
Let $\delta_{x,y,s}$ be the vector
$$\delta_{x,y,s}(j,k,t)=
\begin{cases}
1 & \textup{if $(j,k,t)=(x,y,s)$, and}\\
0 & \textup{otherwise.}
\end{cases}
$$
We have $$S^{-1}(\delta_{x,y,s})(j,k,t)=
\begin{cases}
\frac{1}{n^2} e^{-2\pi ijx/n}e^{-2\pi iky/n} & \textup{if $t=s$, and}\\
0 & \textup{otherwise.}
\end{cases}
$$

{}From (\ref{SAS}) we therefore find, for example,
\begin{equation}\label{Ainverse}
A_1^{-1}((0,0,1),(x,y,3))=\frac{1}{n^2}\sum_{j=0}^{n-1}\sum_{k=0}^{n-1}
\frac{e^{-2\pi i(jx+ky)/n}(b+ae^{-2\pi ij/n})}{
(a+be^{2\pi ij/n})^2e^{-2\pi ij/n}+(c+de^{2\pi ik/n})^2e^{-2\pi ik/n}}
\end{equation}
and
$$A_1^{-1}((0,0,1),(x,y,4))=\frac{1}{n^2}\sum_{j=0}^{n-1}\sum_{k=0}^{n-1}
\frac{e^{-2\pi i(jx+ky)/n}(-i)(d+ce^{-2\pi ik/n})}{
(a+be^{2\pi ij/n})^2z^{-2\pi ij/n}+(c+de^{2\pi ik/n})^2e^{-2\pi ik/n}}.$$
We also have
$A_1^{-1}((0,0,1),(x,y,t))=0$ when $t=1$ or $t=2$.
Similar expressions hold for inverses of $A_2,A_3,A_4$.

An argument identical to the proof of Proposition~\ref{pamn}
(the only difference is the factor $z^{-(jx+ky)}$, which has modulus $1$)
shows that
the parts of the sums (\ref{Ainverse})
over a $\delta$-neighborhood $U_\delta$ of the singularities are $O(\delta)$.

We will show that for all $y$ with
$(1-\varepsilon)n>y>\varepsilon n$ (later we will set $\varepsilon=n^{-1/4}$)
the value $A_1^{-1}((0,0,1),(x,y,3))$ tends to zero as $n\to\infty$ in $W$.
Similar results hold for $A_2$, $A_3$, and $A_4$.
For simplicity of notation let $0,v$
denote the vertices $(0,0,1)$ and $(x,y,3)$.
The equation (\ref{Ainverse}) has the form
$$A_1^{-1}(0,v)=\frac1{n^2}\sum_{j,k}e^{-2\pi i(jx+ky)/n}G_1(j/n,k/n),$$
where $G_1$ is a smooth function on the complement of the region $U_\delta$.
We already know that the sum over $U_\delta$ is $O(\delta)$, so let us replace
$G_1$ by a new function $G_2$
which agrees with $G_1$ outside $U_\delta$ and is zero on $U_\delta$.
We sum by parts over the variable $k$ to get
\begin{eqnarray*}
A_1^{-1}(0,v)=\frac1{n^2}\sum_{j=0}^{n-1}\sum_{k=0}^{n-1}\left[\left(
\sum_{\ell=0}^k e^{-2\pi i\ell y/n}\right)
\left(G_2\left(\frac{j}n,\frac{k}n\right)-
G_2\left(\frac{j}n,\frac{k+1}n\right)\right)\right]\\
\phantom{}+\frac1{n^2}\sum_{j=0}^{n-1}\left(
\sum_{\ell=0}^{n-1}e^{-2\pi i\ell y/n}\right)G_2\left(
\frac jn,0\right) + O(\delta).
\end{eqnarray*}
Since $(1-\varepsilon)n>y>\varepsilon n$,
the sum over $\ell$ of the exponentials is
$$\frac{1-e^{-2\pi i(k+1)y/n}}{1-e^{-2\pi iy/n}}=O\left(\frac1\varepsilon
\right)$$
for each $k$.
The difference
$|G_2({j}/n,{k}/n)-G_2({j}/n,({k+1})/n)|$
is bounded by $1/n$ times
the supremum of $|\p G_1/ \p y|$ on the complement
of $U_\delta$, except that at the points adjacent to the
boundary of $U_\delta$,
the difference is bounded by the supremum of $|G_1|$ near the boundary.

One can check (see the proof of Proposition~\ref{pamn}) that
the sup of $|{\p}G_1/{\p y}|$
on the complement of $U_\delta$ is $O(\delta^{-2})$,
and the supremum of $|G_1|$ on the boundary of $U_\delta$ is
$O(\delta^{-1})$.
Only $O(\delta n)$ pairs $(j,k)$ correspond to points adjacent to
the boundary of $U_\delta$, so we have
\begin{eqnarray*}
|A_1^{-1}(0,v)| &\leq&
\frac1{n^2}\left(\sum_j\sum_kO(\varepsilon^{-1})
O(\delta^{-2})\frac1n\right) +
\frac1{n^2}O(n\delta)O(\delta^{-1})O(\varepsilon^{-1})\\
&&\phantom{}+\frac1{n^2} \sum_j
O(\varepsilon^{-1})O(\delta^{-1})+O(\delta)\\
&=&O\left(\frac1{\varepsilon\delta^2n}\right)+O\left(\frac1{n\varepsilon}
\right)+O\left(\frac1{n\varepsilon\delta}\right)+O(\delta).
\end{eqnarray*}
Choosing $\delta=\varepsilon=n^{-1/4}$, we have
$A^{-1}(0,v)=O(n^{-1/4})$.

A similar argument holds in the case where both $a=b$ and $c=d$.
{}From (\ref{var}) we have
$$\sigma^2(N_a)\leq n^2O(1)+ n^4 o(1)+
\sum_{\textup{edges }k\neq\ell}\sum_{i=1}^4
|A_i^{-1}(v_k,w_\ell)A_i^{-1}(v_\ell,w_k)|,$$
but $A_i^{-1}(v_k,w_\ell)$ and $A_i^{-1}(v_\ell,w_k)$
are $O(n^{-1/4})$ as soon as the edges $k$ and $\ell$
are separated by at least $\epsilon n=n^{3/4}$ in their $y$-coordinates.
Therefore (using $\varepsilon = n^{-1/4}$)
$$\sigma^2(N_a)=O(n^2)+o(n^4) + O(n^2\cdot \varepsilon n^2) + O(n^4 n^{-1/2})=o(n^4).$$
Here the term $O(n^2\cdot\varepsilon n^2)$ comes from
edges whose $y$ coordinate is less than $\varepsilon n$ or greater
than $(1-\varepsilon)n$, and the term $O(n^4 n^{-1/2})$
consists of the remaining pairs of edges.
As a consequence $\sigma^2(N_a)=o(n^4)$.
\end{proof}

One can show from this argument
(using the result of \cite{Kenyon1})
that for $n\in W$ the measures $\mu_{n}$ converge weakly to
a measure $\mu$ on the set of matchings on $\Z^2$.
The result is as follows.

Define
$$
P(2x+1,2y) = \frac{1}{4\pi^2}\int_0^{2\pi}\int_0^{2\pi}
\frac{e^{-i(x\theta+y\phi)}(b+a e^{-i\theta}) \ d\phi \ d\theta}
{(a+be^{i\theta})^2e^{-i\theta}+(c+d e^{i\phi})^2e^{-i\phi}}
$$
and
$$
P(2x,2y+1) = \frac{1}{4\pi^2}\int_0^{2\pi}\int_0^{2\pi}
\frac{e^{-i(x\theta+y\phi)}(-i)(d+c e^{-i\phi}) \ d\phi \ d\theta}
{(a+be^{i\theta})^2e^{-i\theta}+(c+d e^{i\phi})^2e^{-i\phi}}.
$$
Also,
define a {\it colored configuration\/} of dominos as
a configuration of dominos with a checkerboard coloring
of the underlying square grid.

\begin{prop}
\label{couplingfn}
As $n \rightarrow \infty$ within $W$, the probability of finding a
certain colored configuration of dominos in an $(a,b,c,d)$-weighted $n
\times n$ torus converges to $|w \det M|$, where $w$ is the product of
the weights of those dominos and $M_{i,j} = P(v_{i,j})$, with
$v_{i,j} \in \Z^2$ the displacement from the $i$-th white square to
the $j$-th black one.
\end{prop}

Ben Wieland has pointed out that when $ab=cd$, one can write this more
simply.  Define $P'(2x+1,2y) = c(a/b)^x(c/d)^yP(2x+1,2y)$ and
$P'(2x,2y+1) = c(a/b)^x(c/d)^yP(2x,2y+1)$.  One can check that if
$ab=cd$, then the probability we seek is simply the determinant of the
matrix with entries $P'(v_{i,j})$, with no need to multiply by the
products of the weights of the included dominos.  This formulation
makes the conditional uniformity immediately apparent.

\section{The Entropy}
\label{entropysection}

\subsection{Entropy as a function of edge-inclusion probabilities}
Since the set of matchings
on $G_{n}$ is finite, the entropy of a measure $\mu$
on the set of matchings is simply
$$H(\mu)=\sum_{M} -\mu(M)\log\mu(M),$$
where the sum is over all matchings $M$ and $\mu(M)$ is the probability of
$M$ occurring for the measure $\mu$. The {\it entropy per dimer\/} is
by definition
$$\ent(\mu)=\frac1{2n^2}H(\mu)$$
(recall that a matching of $G_{n}$ has $2n^2$ matched edges).

Recall that for real $z$, $L(z)$ is the Lobachevsky function,
defined by \eqref{lobdef}.

\begin{prop}\label{ent1} As $n\to\infty$ in $W$, $\ent(\mu_n)$ converges to
\begin{equation}\label{volume}
\ent(a,b,c,d)=\frac1{\pi}
\left(L(\pi p_a)+L(\pi p_b)+L(\pi p_c)+L(\pi p_d)\right),
\end{equation}
where $p_a,p_b,p_c,p_d$ are given by (\ref{paeqn}).
\end{prop}

\begin{proof}
Let $C(N_a,N_b,N_c,N_d)$ denote
the coefficient of
$$
a^{N_a}b^{N_b}c^{N_c}d^{N_d}
$$
in $Z_{n}$.

As we computed earlier, on the toroidal graph $G_{n}$
the $\mu_{n}$-probability
of an $a$-edge (resp.\ $b$-, $c$-, $d$-edge)
is given by $p_a(n)$ (resp.\ $p_b(n)$,
$p_c(n)$, $p_d(n)$).
The expected number of $a$-edges is $\overline{N}_a\stackrel{\textup{def}}{=}
\E(N_a)= 2n^2p_a(n)$.

Let
\begin{eqnarray*}
U_{\varepsilon} &=& \{(N_a,N_b,N_c,N_d):
|N_a-\overline{N}_a|<\varepsilon
\overline{N}_a,|N_b-\overline{N}_b|<\varepsilon \overline{N}_b,\\
&& |N_c-\overline{N}_c|<\varepsilon \overline{N}_c,|N_d-
\overline{N}_d|<\varepsilon \overline{N}_d\}.
\end{eqnarray*}
Let $V_{\varepsilon}$ be the corresponding set of matchings, i.e., those
where the corresponding quadruples $(N_a,N_b,N_c,N_d)$ are in $U_\varepsilon$.
Because the variance is $o(n^4)$,
for all $\varepsilon,\varepsilon'>0$ there exists $n_0$ such that
for $n\in W$ greater than $n_0$ we have
$$\sum_{U_\varepsilon}
C(N_a,N_b,N_c,N_d)a^{N_a}b^{N_b}c^{N_c}d^{N_d}
\geq (1-\varepsilon')Z_{n}(a,b,c,d).$$

Note that if $p_j$ are probabilities and $p_1+\dots+p_k<\varepsilon'<1$,
then $$-\sum p_j\log p_j\leq
-\varepsilon'\log\left(\frac{\varepsilon'}{k}\right)$$
(since the left-hand side is maximized when the $p_j$'s are equal).
Thus for the entropy we may write
\begin{eqnarray*}
H(\mu)& =& -\sum_{M\not\in V_\varepsilon} \mu(M)\log\mu(M)
-\sum_{M\in V_\varepsilon} \mu(M)\log\mu(M)\\
&=& -\sum_{M\not\in V_\varepsilon} \mu(M)\log\mu(M)
-\sum_{M\in V_\varepsilon}
\mu(M)\log\left( \frac{a^{N_a}b^{N_b}c^{N_c}d^{N_d}}{Z_{n}}\right)\\
&=&
O\left(\varepsilon'\log\left(\frac{\varepsilon'}
{\textup{constant}^{n^2}}\right)\right)\\
&&\phantom{}-\sum_{M\in V_\varepsilon}\mu(M)\log\left(
\frac{a^{\overline{N}_a}
b^{\overline{N}_b}c^{\overline{N}_c}d^{\overline{N}_d}}{Z_{n}}
a^{N_a-\overline{N}_a}
\dots d^{N_d-\overline{N}_d}\right),
\end{eqnarray*}
but for $M\in V_\varepsilon$, we have
$\log(a^{N_a-\overline{N}_a})<\varepsilon \log(a^{\overline{N}_a})$,
and similarly for $b,c,d$, so
$$H(\mu)=
\left[\sum_{M\in V_\varepsilon}\mu(M)\right]
\left(-\log(a^{\overline{N}_a}
b^{\overline{N}_b}c^{\overline{N}_c}d^{\overline{N}_d})
(1+O(\varepsilon))+ \log Z_{n}\right)
+O(n^2\varepsilon'\log\varepsilon').$$

Note also that  $\sum_{M\in V_\varepsilon}\mu(M)\geq 1-\varepsilon'$.
Letting $\varepsilon,\varepsilon'\to 0$ as $n\to\infty$ we have finally that
the limiting entropy per dimer is
\begin{eqnarray*}
\ent(a,b,c,d) &=& \lim_{n\to\infty}\frac{1}{2n^2}\left(\log Z_{n}(a,b,c,d)-
\log(a^{\overline{N}_a}
b^{\overline{N}_b}c^{\overline{N}_c}d^{\overline{N}_d})\right)\\
&=&\log Z(a,b,c,d)- p_a \log(a)-p_b\log(b)-p_c\log(c)-p_d\log(d).
\end{eqnarray*}

Without loss of generality we may assume that $a\geq b$ and $c\geq d$;
then from (\ref{thetapcpd}) we
have $\theta_0=\pi-\pi(p_c-p_d)$.
Plugging in from (\ref{Zformula}) now gives
\begin{equation}\label{entent}
\ent(a,b,c,d)=\frac{1}{4\pi}\int_{-\theta_0}^{\theta_0}\log|\beta(z)|\:d\theta+
\frac{1-p_c-p_d}{2}\log(cd)-p_a\log(a)-p_b\log(b).
\end{equation}

To prove the equivalence of this formula and (\ref{volume}),
we show that they agree when $a=b=c=d$, and show that
their partial derivatives are equal for all $a,b,c,d$.

Formula (\ref{entent}) gives
\begin{eqnarray*}
\ent(1,1,1,1) &=& \frac1{4\pi}\int_{-\pi}^{\pi}
\log\left(2+\cos{\theta}+\sqrt{(2+\cos{\theta})^2-1}\right)d\theta\\
&=&\frac1{4\pi}\int_{-\pi}^{\pi}2\log\left(\cos\left(\frac{\theta}2\right)+
\sqrt{\cos^2\left(\frac{\theta}2\right)+1}\right)d\theta.
\end{eqnarray*}
This is two times the value of the entropy per site
given in formula (17) of \cite{Kast1}, as it should be since the entropy
$\ent(1,1,1,1)$ as we defined it is the entropy per dimer.
Kasteleyn also shows that this value equals $2G/{\pi}$, where
$G$ is Catalan's constant
$$G=1-\frac1{3^2}+\frac1{5^2}-\frac1{7^2}+\cdots.$$
{}From the expansion
$$L(x)= \frac12\sum_{k=1}^\infty \frac{\sin(2kx)}{k^2}$$
(see \cite{Mil}) we have $2G/\pi =(4/\pi)L(\pi/4)$,
so the two formulas agree when $a=b=c=d$.

It remains to compute the derivatives.  Equation \eqref{volume} is
symmetric under the full symmetry group $S_4$, and \eqref{entent} is
by definition symmetric under the operations of exchanging under the
operations of exchanging $a$ and $b$, exchanging $c$ and $d$, and
exchanging $a,b$ with $c,d$.  These operations are transitive on $\{a,
b, c, d\}$, so
it suffices to show equality of the derivatives with respect to $a$.
We have
\begin{eqnarray*}
\frac{\partial}{\partial a}\ent(a,b,c,d) &=&
\frac{\partial}{\partial a}(\log
Z-p_a\log(a)-p_b\log(b)-p_c\log(c)-p_d\log(d))\\
&=& -\log(a)\frac{\partial p_a}{\partial a}
-\log(b)\frac{\partial p_b}{\partial a}-
\log(c)\frac{\partial p_c}{\partial a}-\log(d)\frac{\partial
p_d}{\partial a}
\end{eqnarray*}
(recall that $p_a=({a}/{Z}){\partial Z}/{\partial a}$).

On the other hand when $x$ is one of $a,b,c,d$ we have
$$\frac{\partial}{\partial a} \frac{1}{\pi}L(\pi p_x)=-\log(2\sin(\pi p_x))
\frac{\partial p_x}{\partial a}.$$
Taking $\p/\p a$ of (\ref{volume}) and recalling that
$a=\sin(\pi p_a)$, etc., gives
$$-\log(2a)\frac{\partial p_a}{\partial a}
-\log(2b)\frac{\partial p_b}{\partial a}-
\log(2c)\frac{\partial p_c}{\partial a}-
\log(2d)\frac{\partial p_d}{\partial a}.$$
Since
$$\log(2)\frac{\partial}{\partial a}(p_a+p_b+p_c+p_d) = 0,$$
the proof is complete.
\end{proof}

As explained in the proof of Theorem~\ref{precise},
the entropy converges for all $n$,
not just for $n\in W$:
\begin{thm} \label{ent}
As $n\to\infty$
the entropy per edge of matchings on $G_{n}$
of the measure $\mu_{n}(a,b,c,d)$ converges to
$$\ent(a,b,c,d)=\frac1{\pi}
\left(L(\pi p_a)+L(\pi p_b)+L(\pi p_c)+L(\pi p_d)\right),$$
where $p_a,p_b,p_c,p_d$ are given by (\ref{paeqn}).
\end{thm}

Intriguingly, this formula can be used to show
that when each of $a,b,c,d$
is less than the sum of the others (the only case in which the entropy
$\ent(s,t)$ is non-zero),
$\ent(s,t)$ is equal to $1/\pi$ times the volume
of a three-dimensional ideal hyperbolic pyramid,
whose vertices in the upper-half-space model are
the vertex at infinity and
the four vertices of the cyclic quadrilateral
of Euclidean edge lengths $a,c,b,d$ in cyclic order
(otherwise the entropy is $0$).
We have no conceptual explanation for this coincidence.

We do not even fully understand why
the limiting behavior of the measures $\mu_n$,
viewed as a function of $a,b,c,d$,
turns out to be symmetrical in its four arguments.
This symmetry is not merely combinatorial,
since it emerges only in the limit
as the size of the torus goes to infinity.
R.~Baxter suggests (in personal communication)
that it is almost certainly the same symmetry that occurs in the
checkerboard Ising model.  In that setting it can be proved using the
Yang-Baxter relation (see \cite{JM,MR}).

\subsection{Entropy as a function of tilt}
\label{entasfnofslope}
Given a tilt $(s,t)$ satisfying $|s|+|t|<2$,
we claim that
there is a unique (up to scale) $4$-tuple
of weights $a,b,c,d$ satisfying the conditional uniformity
property $ab=cd$ and such
that the average tilt of $\mu_{a,b,c,d}$ is $(s,t)$.
To determine $a,b,c,d$, we note that  $p_a,p_b,p_c,p_d$
are determined by the equations (\ref{fx})--(\ref{sinsin}).
To solve these equations, note that (\ref{sinsin}) can be written
$$\cos(\pi(p_a-p_b))-\cos(\pi(p_a+p_b))=
\cos(\pi(p_c-p_d))-\cos(\pi(p_c+p_d)),$$
and so
$$\cos(\pi t/2)-\cos(\pi(p_a+p_b))=\cos(\pi s/2)-\cos(\pi-\pi(p_a+p_b)),$$
giving
$$2\cos(\pi(p_a+p_b))=\cos(\pi t/2)-\cos(\pi s/2).$$
This combined with (\ref{s}) gives (\ref{probdef}),
where the values of $\cos^{-1}$ are taken from $[0,\pi]$
(to see why, notice that $\cos^{-1}((\cos(\pi t/2) - \cos(\pi s/2))/2) =
\pi (p_a+p_b)$, which is between $0$ and $\pi$).

Finally we can determine $a,b,c,d$ as functions of the tilt by
$$a=\sin(\pi p_a),~~b=\sin(\pi p_b),~~c=\sin(\pi p_c),~~ d=\sin(\pi p_d).$$

\section{Concavity of the Entropy}
\label{concavesection}
\begin{thm}
\label{concavethm}
The entropy per edge $\ent(s,t)$
is a strictly concave function of $s,t$ over the range
$|s|+|t|\leq 2$.
\end{thm}

\begin{proof}
We show that the Hessian (matrix of second derivatives) is negative definite,
that is, $\ent_{ss}(s,t)<0,$
$\ent_{tt}(s,t)<0$, and $\ent_{ss}(s,t)\ent_{tt}(s,t)-\ent_{st}(s,t)^2>0$,
for all $s,t$, except at the four points
$(s,t)=(\pm2,0)$ or $(0,\pm2)$.

A computation using  Theorem~\ref{ent} and equations (\ref{Peqns}) gives
\begin{eqnarray*}
\ent_t(s,t)&=&\frac{\p\ent}{\p p_a}\frac{\p p_a}{\p t}+\dots+
\frac{\p\ent}{\p p_d}\frac{\p p_d}{\p t}\\
&=&-\frac14\log\left(\frac{\sin(\pi p_a)}{\sin(\pi p_b)}\right).
\end{eqnarray*}
A second differentiation yields
\begin{eqnarray*}
\frac{\partial^2 \ent(s,t)}{\partial t^2} &=&
\frac{-\pi}{32\sin(\pi(p_a+p_b))\sin(\pi p_a)\sin(\pi p_b)} \times\\
&&
\left(\sin^2\left(\frac{\pi s}{2}\right)+\frac{(\cos(\frac{\pi
s}2)+\cos(\frac{\pi
t}2))^2}{2}\right),
\end{eqnarray*}
and this quantity is strictly negative except at the points
$(s,t)=(\pm2,0),(0,\pm 2)$.

A similar calculation holds for $\ent_{ss}(s,t)$:
\begin{eqnarray*}
\frac{\partial^2 \ent(s,t)}{\partial s^2} &=&
\frac{-\pi}{32\sin(\pi(p_a+p_b))\sin(\pi p_a)\sin(\pi p_b)} \times\\
&&
\left(\sin^2\left(\frac{\pi t}{2}\right)+\frac{(\cos(\frac{\pi
s}2)+\cos(\frac{\pi
t}2))^2}{2}\right).
\end{eqnarray*}

We have
$$\ent_{st}(s,t)=\frac{-\pi}{32}\frac{\sin(\pi s/2)\sin(\pi t/2)}{\sin(\pi p_a)
\sin(\pi p_b)\sin(\pi(p_a+p_b))}.$$
Finally,
$$\ent_{ss}\ent_{tt}-\ent_{st}^2=
\left
(\frac{\pi}{32\sin(\pi(p_a+p_b))\sin(\pi p_a)\sin(\pi p_b)}
\right)^2\times$$
$$
\left[\left(\sin^2\left(\frac{\pi t}{2}\right)+
\frac{(\cos(\frac{\pi s}2)+\cos(\frac{\pi t}2))^2}{2}\right)
\left(\sin^2\left(\frac{\pi s}{2}\right)+\frac{(\cos(\frac{\pi s}2)+\cos(
\frac{\pi t}2))^2}{2}\right)\right.$$
$$
\hspace{2in} \left.
-\sin^2\left(\frac{\pi s}2\right)\sin^2\left(\frac{\pi t}2\right)\right],$$
which is clearly positive.
\end{proof}

\section{Proof of Theorem~\protect{\ref{pa}}}
\label{paproof}
In what follows, we must be careful to distinguish
the differential $dw$
from the product $d\cdot w$
(we will write the product as $wd$ to avoid confusion).
Let $z=e^{i\theta}$, $w=e^{i\phi}$
and $r(z)=cd+{(a+bz)^2}/({2z})$ as before. Then (see (\ref{q}))
\begin{eqnarray*}
p_a&=&\frac{a}{4\pi^2}\int_{S^1}\int_{S^1}\frac{w(b+a/z)}{
c^2+2rw+w^2d^2}\cdot\frac{dw}{iw}\cdot\frac{dz}{iz}\\
&=&\frac{-a}{4\pi^2}\int_{S^1}\frac{(b+a/z)dz}{z}
\int_{S^1}\frac{dw}{(wd-\alpha c)(
wd-\beta c)}.
\end{eqnarray*}
Recall (see the second-to-last paragraph of Subsection \ref{roots})
that when $|z|=1$, $|\beta(z)c|>d$ always,
and $|\alpha(z) c|<d$ if
and only if
$\theta\in(-\theta_0,\theta_0)$
(remember that we defined $\theta_0=\pi$ in the case $a\geq b+c+d$).
If $|\alpha(z) c|<d$, then the residue of $((wd-\alpha c)(wd-\beta c))^{-1}$
is $$\frac{1}{(\alpha-\beta)cd}.$$

Thus we have
\begin{eqnarray*}
p_a &=&\frac{-a}{4\pi^2}\frac{2\pi i}{cd}\int_{\theta=-\theta_0}^{
\theta=\theta_0}
\frac{b+a/z}{z(\alpha-\beta)}\,dz\\
&=& \frac{-ai}{2\pi cd}\int_{\theta=-\theta_0}^{\theta=\theta_0}
\frac{(bz+a)\,dz}{2z^2\sqrt{(\frac{r}{cd}-1)(\frac{r}{cd}+1)}},
\end{eqnarray*}
and recalling the definition of $r$ and simplifying yields
$$p_a=\frac{-ai}{2\pi}\int_{\theta=-\theta_0}^{\theta=\theta_0}\frac{dz}
{z \sqrt{(a+bz)^2+4zcd}}.$$
We don't have to worry about keeping track of the sign of the square
root since we know that we want $p_a\geq0$.
In fact we only need to be careful about
the sign when we get to (\ref{sign}), below.

This integral can be explicitly evaluated, giving
$$p_a=\frac{i}{2\pi}\left[\log\left(
\frac{a^2+(ab+2cd)z+a\sqrt{(a+bz)^2+4zcd}}{z}\right)
\right]_{e^{-i\theta_0}}^{e^{i\theta_0}}.$$
This expression can be simplified using the variable
(or rather, one of the two variables) $w=w(z)$ such that
${(a+bz)^2}/{z}+{(c+wd)^2}/{w}=0$:
the expression under the square root is then
$$(a+bz)^2+4zcd= z\left(-\frac{(c+wd)^2}{w}+4cd\right)=
\frac{-z}{w}(c-wd)^2.$$
Plugging this in yields
\begin{eqnarray}\nonumber
 p_a&=&\frac{i}{2\pi}\left[\log\left(\frac{a^2}{z}+ab+2cd+a(c-wd)
\frac{i}{\sqrt{wz}}\right)\right]_{e^{-i\theta_0}}^{e^{i\theta_0}}\\
\nonumber
&=&\frac{i}{2\pi}\left[\log\left(2cd+\frac{a}{\sqrt{z}}\left\{
\frac{a}{\sqrt{z}}+b\sqrt{
z}+i\left(\frac{c}{\sqrt{w}}+d\sqrt{w}\right)\right\}
-2ida\frac{\sqrt{w}}{\sqrt{z}}
\right)\right]_{e^{-i\theta_0}}^{e^{i\theta_0}}.\\
\label{sign}
\end{eqnarray}
Up until (\ref{sign}), changing the sign of the square root
will only change the sign of the integral.
In (\ref{sign}), we choose the sign of ${\sqrt{w/z}}$
(or, what is the same, the sign of $\sqrt{wz}$)
so that the expression in curly brackets is zero
(cf.~(\ref{mod1})). We then have
$$p_a=\frac{i}{2\pi}
\left[\log\left(-2d\left(-c+ia\sqrt{\frac{w}{z}}\right)\right)
\right]_{e^{-i\theta_0}}^{e^{i\theta_0}}.
$$
Had we chosen the other sign we would have gotten a similar expression
with $c$ and $d$ interchanged.

\begin{figure}[htbp]
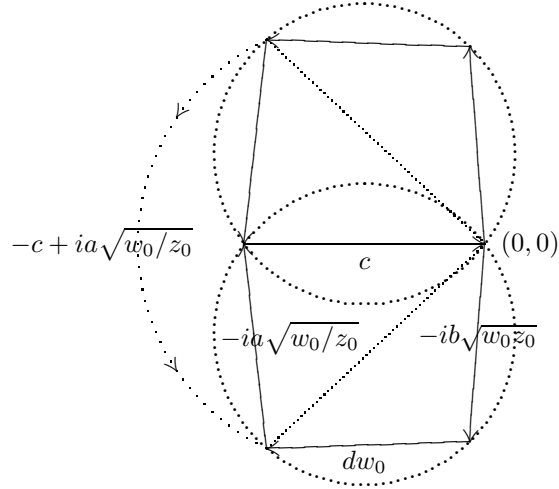

\begin{center}
\input integ.tex
\end{center}
\caption{\label{2quads}The integral defining $p_a$. Note that the
angle between the two dotted rays at the origin is $\theta_a$. (Recall
that $\theta_a=2\pi p_a$ is
the angle of arc cut off by the edge $a$.)}
\end{figure}

But now Figure \ref{2quads} (whose lower quadrilateral is a rotation
of Figure \ref{cycquad}, and whose upper quadrilateral is the reflection
of the lower across the edge $c$) shows that, as $\theta$ runs from
$-\theta_0$ to $\theta_0$, the quantity
$-c+ia\sqrt{w/z}$ sweeps out an angle of $\theta_a$, the angle
of arc cut off by edge $a$ in a cyclic quadrilateral of edge lengths $c,a,d,b$.
Thus $p_a={\theta_a}/({2\pi})$.
In the case $a>b+c+d$, one can similarly show that $-c+ia\sqrt{w/z}$
sweeps out an angle of $2\pi$, and thus $p_a=1$.

Finally, to prove the formula (\ref{paeqn}),
recall the well-known formula for the radius $r$ of the circumcircle
of a triangle of sides $a,b,c$:
\begin{equation}\label{circum}
r^2=\frac{a^2b^2c^2}{(a+b+c)(a+b-c)(a-b+c)(-a+b+c)}.
\end{equation}
{}From this one can compute, for a cyclic
quadrilateral of sides $a,c,b,d$, the length $s$ of the
diagonal having the $b$ and $c$ edges on the same side:
$$s^2=\frac{(ab+cd)(ac+bd)}{ad+bc}.$$
Plugging this value back into the formula (\ref{circum}) with $a,b,c$
replaced with $a,d,s$ gives
$$r^2=\frac{(ab+cd)(ac+bd)(ad+bc)}{(a+b+c-d)(a+b-c+d)(a-b+c+d)(-a+b+c+d)},$$
from which (\ref{paeqn}) follows by
$p_a={\pi}^{-1}\textup{sin}^{-1} ({a}/({2r}))$.

\section{The PDE}
\label{PDEsection}
Under the assumption that the entropy-maximizing function $f$ is
$C^2$, the Euler-Lagrange equation for $f$ is
$$\frac{d}{dx}(\ent_s(f_x,f_y))+\frac{d}{dy}(\ent_t(f_x,f_y))=0,$$
where $\ent_s,\ent_t$ are the partial derivatives with respect to the first
and second variable, respectively.
This equation holds only at points where the tilt $(f_x,f_y)$ is
non-extremal, i.e.,
satisfies $|f_x|+|f_y| < 2$;  otherwise, perturbing $f$ will mess up
the Lipschitz condition.
In case $f$ is only $C^1$,
this equation still holds in a distributional sense:
it is true when integrated against any smooth test function $g$
vanishing on the boundary (and such that $f+\varepsilon g$
is 2-Lipschitz for sufficiently small $\varepsilon>0$).
In such a case $f$ is called a weak solution \cite{GT}.

We computed in the proof of Theorem~\ref{concavethm} that
\begin{eqnarray}\label{Epq}
\ent_s(s,t)&=&-\frac14\log\left(\frac{\sin(\pi p_d)}{\sin(\pi p_c)}\right),
\textup{ and}
\nonumber\\
\ent_t(s,t)&=&-\frac14\log\left(\frac{\sin(\pi p_a)}{\sin(\pi p_b)}\right).
\end{eqnarray}

Plugging in from (\ref{Peqns}) and simplifying yields
the following PDE:

\begin{thm}
\label{pdethm}
At the points where the entropy-maximizing function $f$
is $C^2$ and has non-extremal tilt, it satisfies the PDE
\begin{eqnarray*}
\left(2\left(1-D^2\right)-\sin^2\left(\frac{\pi
f_x}2\right)\right)f_{xx}
&+& 2\sin\left(\frac{\pi f_x}2\right)
\sin\left(\frac{\pi f_y}2\right)f_{xy}\\
&+& \left(2\left(1-D^2\right)-\sin^2\left(\frac{\pi
f_y}2\right)\right)f_{yy}=0,
\end{eqnarray*}
where $D=\frac12\left(\cos\left(\pi f_x/2\right)-\cos\left(\pi
f_y/2\right)\right)$.
\end{thm}

\section{Conjectures and Open Problems}
\label{conjs}

In this article, $p_a$, $p_b$, $p_c$, and $p_d$ were defined in
relation to the dimer model on an $n \times n$ torus, in the
thermodynamic limit as $n \rightarrow \infty$.  We have proved no
corresponding interpretation of these quantities for the
thermodynamic limit of planar regions.  However, our results
on asymptotic height function associated with large regions
imply that, in a patch of a large region where the associated
asymptotic height function $f$ satisfies
$({\partial f}/{\partial x},{\partial f}/{\partial y})=(s,t)$,
the local density of $a$-edges minus the local density of $b$-edges
equals $p_a-p_b$, and likewise for the $c$- and $d$-edges.
To be precise here, one would define local densities as averages
over mesoscopic patches within the tiling (which we recall are
defined as patches whose absolute size goes to infinity but
whose relative size goes to zero).
\begin{conj}
\label{denseconj}
The local densities of $a$-edges, $b$-edges, $c$-edges, and $d$-edges
are given by $p_a$, $p_b$, $p_c$, and $p_d$, respectively, in the
thermodynamic limit.
\end{conj}
\noindent
Here our use of the phrase ``thermodynamic limit'' carries along with
it the supposition that we are dealing with an infinite sequence of
ever-larger
planar
regions whose normalized boundary height functions converge
to some particular boundary asymptotic height function, and that the
subregions we are studying stay away from the boundary by at least some
mesoscopic distance.

The local density of $a$-edges is just the average of the
inclusion probabilities of all the $a$-edges within a mesoscopic patch.
The authors have empirically observed that such averages do not arise
from the smoothing out of genuine fluctuations.  Rather, it seems that
apart from degenerate cases, all the $a$-edges within a patch have
roughly the same inclusion probability.  These degenerate cases occur
in patches where the asymptotic height function is not smooth (so that
$(s,t) = ({\partial f}/{\partial x}, {\partial f}/{\partial y})$
is undefined), or where the tilt is nearly extremal (i.e., $|s|+|t|$ is
close to 2).
Hence we believe:
\begin{conj}
In the thermodynamic limit, the probability of seeing a domino
in a particular location is given by the suitable member of the
4-tuple $(p_a,p_b,p_c,p_d)$,
wherever the tilt $(s,t)$ (given by the partial derivatives of
the entropy-maximizing height function) is defined and satisfies
$|s|+|t|<2$.
\label{probconj}
\end{conj}

This is the conjectural interpretation of $p_a,p_b,p_c,p_d$ that was
alluded to in Subsection~\ref{statementsec}.  We are two removes from
being able to prove it in the sense that we do not even know how to
prove Conjecture~\ref{denseconj}.

The restriction $|s|+|t|<2$ deserves some comment.
It is not hard to devise a large region
composed of long diagonal ``herringbones''
that has only one tiling (see the final section of \cite{CEP}).
For such regions, the edge-inclusion probabilities can be made to
fluctuate erratically between 0 and 1.
Such regions have asymptotic height functions
in which the tilt is extremal everywhere it is defined,
so we can rule out such behavior on the basis of tilt.

Incidentally,
a conjecture analogous to Conjecture~\ref{probconj} can be made for
the case of lozenge tilings (which can be studied using the methods of
this paper, by setting one edge weight equal to $0$).
Here the analogue of Conjecture~\ref{denseconj} is
actually a theorem; that is, for lozenges there are only three kinds
of orientations of tiles, so that their relative frequencies, which
jointly have two degrees of freedom, both determine and are determined
by the local tilt $(s,t)$ of the asymptotic height function.

Straying further into the unknown, we might inquire about the
probabilities of finite (colored) configurations of tiles.
Here again we are guided by the {\it Ansatz\/} that what is
true for large tori should be true for large finite regions as
well.  Given any tilt $(s,t)$ satisfying $|s|+|t| \leq 2$, choose
weights $a,b,c,d$ satisfying $ab=cd$ and giving tilt $(s,t)$
in accordance with our earlier formulas, and use
the formula in Proposition~\ref{couplingfn} to define a measure
$\mu_{s,t}$ on the space of domino tilings of the plane.  This
measure is invariant under color-preserving translations and
satisfies the property of ``conditional uniformity''---given
any finite region, the conditional distribution upon fixing a
tiling of the rest of the plane is uniform.   Proposition~\ref{couplingfn}
invites us to surmise:
\begin{conj}
\label{torusconj}
Let $|s|+|t| \leq 2$, and let $a,b,c,d$ be weights satisfying $ab=cd$
such that the average height function for weighted torus tilings
has tilt $(s,t)$.  Then for any colored configuration of dominos, the
probability of finding it in a specified
location in a random $n \times n$
torus tiling converges as $n \rightarrow \infty$
to the value given by the measure $\mu_{s,t}$.
\end{conj}
Proposition~\ref{couplingfn} approaches this claim,
but it restricts $n$ to lie in a large subset $W$ of the integers.

The measures $\mu_{s,t}$ have positive entropy whenever $|s|+|t| <
2$.  Furthermore, the coupling function calculations in the proof of
Proposition~\ref{sigmasquared}  show that these measures are mixing
(correlations between distant cylinder sets tend to $0$) and hence
ergodic.  We believe that these measures can be characterized uniquely
by the properties mentioned so far, although we cannot prove it.

\begin{conj}
Every ergodic, conditionally uniform measure on the set of tilings of
the plane that is invariant under color-preserving translations and
has positive entropy is of the form $\mu_{s,t}$ for some $(s,t)$
satisfying $|s|+|t|<2$.
\end{conj}

Assuming the truth of Conjecture~\ref{torusconj},
it would be natural to advance a further claim that would
come close to being the final, definitive answer to
the question Kasteleyn raised nearly four decades ago:
\begin{conj}
\label{lastconj}
In the thermodynamic limit for a sequence of finite regions converging
to a fixed shape,
the probability of seeing any colored
configuration of dominos is given by the measure $\mu_{s,t}$,
wherever the tilt $(s,t)$
given by the variational principle
is defined and satisfies
$|s|+|t|<2$.
\end{conj}

Finally, we discuss the seemingly miraculous geometric interpretation
of our formula for the entropy.
As was pointed out in Section~\ref{introsec}, when each of $a,b,c,d$ is
less than the sums of the others,
the (asymptotic) entropy of
torus tilings with weights $a,b,c,d$ equals $1/\pi$ times the volume
of a three-dimensional ideal hyperbolic pyramid,
whose vertices in the upper-half-space model are
the vertex at infinity and
the four vertices of the cyclic quadrilateral
of Euclidean edge lengths $a,c,b,d$ (in cyclic order).
Notice that we needn't assume any particular scaling for the weights,
because homotheties $(x,y,z) \mapsto (rx,ry,rz)$
are isometries of hyperbolic space.

\begin{prob}
What do tilings have to do with hyperbolic geometry?
\end{prob}
\noindent
Explaining that connection is one of the most intriguing open problems
in this area.

\section*{Acknowledgements}

We thank Matthew Blum, Ben Raphael, Ben Wieland, and David Wilson for
helpful discussions.

\end{document}